\documentclass[11pt]{article}
\usepackage{amsmath,amssymb,latexsym,amsthm}

\usepackage{amsmath}
\usepackage{color}
\usepackage{graphicx}
\usepackage{graphics}
\usepackage{subfigure}

\usepackage{subfigure}
\usepackage{epstopdf}
\usepackage{array}
\usepackage{tabularx}
\usepackage{multirow}
\usepackage{graphicx}
\usepackage{multirow}
\usepackage{amssymb}
\usepackage{mathrsfs}
\usepackage{amsmath}
\usepackage{enumerate}
\usepackage{multirow}

\def\LARGE{\Large}
\def\Large{\large}

\begin{document}

\begin{center}
{\LARGE{\bf {\Large Cooperative hunting in a discrete
predator-prey system}}}
\end{center}

\begin{center} Yunshyong Chow$^1$,   Sophia R.-J.  Jang$^2$, and
Hua-Ming Wang$^3$ \end{center}

\noindent 1. Institute of Mathematics, Academia Sinica, Taipei
10617, Taiwan

\noindent 2. Department of Mathematics and Statistics, Texas Tech
University, Lubbock, TX 79409, USA

\noindent 3. Department of Statistics, Anhui Normal University,
Wuhu 241003, Anhui, China

\vspace{0.1in}

\noindent\textbf{Abstract.} We propose and investigate a
discrete-time predator-prey system with cooperative hunting in the
predator population. The model is constructed from the classical
Nicholson-Bailey host-parasitoid system with density dependent
growth rate. A sufficient condition based on the model parameters
for which both populations can coexist is derived, namely that the
predator's maximal reproductive number exceeds one. We study
existence of interior steady states and their stability in certain
parameter regimes. It is shown that the system behaves
asymptotically similar to the model with no cooperative hunting if
the degree of cooperation is small. Large cooperative hunting,
however,  may promote persistence of the predator for which the
predator would otherwise go extinct if there were no cooperation.

\vspace{0.1in}


\vspace{0.1in}


\section{Introduction}
Cooperation between individuals of  social animals is frequently
observed and widespread in biological systems. For example,
carnivores such as wolves, wild dogs and lions   often work
together to capture and kill their preys \cite{scheel}. Other
organisms such as spiders, birds and ants also seek and  attack
prey collaboratively \cite{uetz}. However, there are only a few
mathematical models constructed to study such a biological
phenomenon.

Previous research  on cooperative hunting includes Berec
\cite{berec} who uses ordinary differential equations to model
predator-prey interactions with a Holling type II functional
response. Due to this functional response, Berec studies the
effects of cooperative hunting relative to population
oscillations. Cosner et al. \cite{cosner} on the other hand
propose models of partial differential equations to explore the
effects of predator aggregation when predators encounter a cluster
of prey.  Recently, Alves and Hilker \cite{alves} use models of
ordinary differential equations of predator-prey interactions with
cooperative hunting in predators to investigate  impacts of
cooperative hunting. It is concluded that cooperative hunting can
improve persistence of the predator but may also promote a sudden
collapse of the predator. In addition, this research suggests that
cooperative hunting is a mechanism for inducing Allee effects in
predators.

Ever since the pioneer work of May \cite{may}, mathematical models
of difference equations have played  important roles in the
understanding of population interactions. There are many
populations with non-overlapping generations and discrete-time
models are more appropriate to describe such populations.
Additionally, data of ecological studies are usually collected in
discrete formats. Motivated by these, the goal of this study is to
propose and  investigate the effects of cooperative hunting among
predators upon predator-prey interactions in the discrete-time
setting. Our model derivation is built on the well-known
Nicholson-Bailey model with density-dependent prey growth rate.
Based on the stability of the boundary equilibria, we provide a
set of sufficient conditions for population coexistence, where the
conditions do not depend on the cooperative hunting. We show that
the system has the same asymptotic dynamics as the model of no
cooperative hunting if the degree of cooperation is small. If the
degree of cooperative hunting is large, then the system may
support two coexisting steady states for which the predator would
otherwise go extinct if there were no cooperation in this
parameter regime. Numerical simulations are presented to confirm
our analytical findings and to further our understanding of the
predator-prey interactions.

In the following section, model derivation and persistence of the
populations are presented. Section 3 provides results on the
existence and the number of interior steady states. Asymptotic
dynamics and local stability of the interior steady states are
given in Section 4. The final section summarizes results and
provides conclusions.

\section{Model derivation and persistence of populations}
\setcounter{equation}{0} Let $N(t)$ and $P(t)$ denote the hosts
and parasitoids in generation $t$ respectively, $t=0, 1, \cdots$.
In the classical Nicholson-Bailey model \cite{allen}, the number
of encounters between hosts  and parasitoids  in generation $t$ is
assumed to follow the law of  mass action,  $aN(t)P(t)$, where the
constant $a>0$ denotes searching efficiency of the parasitoids. It
is also assumed that the number of encounters is distributed
randomly and follows a Poisson distribution with probability
$p(n)=\cfrac{e^{-\mu}\mu^n}{n!}$, where $n=0, 1, 2, \cdots$ is the
number of encounters and $\mu$ is the average of encounters per
host per generation. It follows that
$\mu=\cfrac{aN(t)P(t)}{N(t)}=aP(t)$ and thus $1-p(0)=1-e^{-aP(t)}$
is the probability of an individual host being parasitized in
generation $t$ since only the first encounter results in
parasitism. The well known Nicholson-Bailey model is given by
\begin{equation}\label{NB}\left\{\begin{array}{ll}
N(t+1)= r N(t)\displaystyle e^{-aP(t)}\\[1ex]
P(t+1)=c N(t)\left(1-e^{-P(t)}\right),
\end{array}\right.
\end{equation}
where all of the parameters are positive constants. Notice that
the host  population grows exponentially in \eqref{NB} and the
unique interior steady state is always unstable when it exists
\cite{allen}.

In the context of predator-prey interactions, we let $x(n)$ and
$y(n)$ denote respectively the prey and predator populations at
time $n=0, 1, 2, \cdots$. In the absence of cooperative hunting
and by applying a similar argument as in the derivation of
Nicholson-Bailey model, the probability of a prey that escaped
from predation at time $n$ is $e^{-ay(n)}$. With cooperative
hunting, the number of encounters between prey and predators at
time $n$ becomes $ax(n)y(n)(1+\alpha y(n))$, where $\alpha\geq 0$
denotes degree of cooperative hunting. There is no cooperation
among  predators if $\alpha=0$ and the cooperation is stronger if
$\alpha$ is larger.  It follows that  the probability of an
individual prey escaped from being preyed upon at time $n$ is
$e^{-ay(n)(1+\alpha y(n))}$. The probability is smaller due to
cooperation among predators.

To avoid the perpetual instability of the interior steady state in
\eqref{NB}, we modify the density-independent growth rate given in
\eqref{NB} by assuming that the per capita growth rate of the prey
is density-dependent and is modeled by the Beverton-Holt equation.
Putting these together, the interaction between prey and predator
populations is described by the following difference equations:
\begin{equation}\label{gen}\left\{\begin{array}{ll}
x(n+1)= \lambda x(n)g_0(x(n))\displaystyle e^{-ay(n)(1+\alpha y(n))}\\[1ex]
y(n+1)=\beta x(n)\left(1-e^{-ay(n)(1+\alpha y(n))}\right)
\end{array}\right.
\end{equation}
with nonnegative initial conditions, where $\lambda
g_0(x)=\lambda/(1+kx)$, $\lambda, k>0$, is the prey's per capita
growth rate. The parameter $\alpha\geq 0$ denotes cooperative
hunting of the predator if $\alpha>0$, and $a>0$ is the searching
efficiency of the predator. Further, $\beta>0$ is the predator
conversion for each prey consumed.

We   nondimensionalize system \eqref{gen} by letting
\begin{equation}\label{dim}\tilde x=kx,\ \tilde y=ay,
\tilde\beta=\cfrac{\beta a}{k}, \mbox{ and } \tilde
\alpha=\cfrac{\alpha}{a}.\end{equation} Ignoring the tildes,
\eqref{gen} is converted into the following system with only three
parameters
\begin{equation}\label{rmodel}\left\{\begin{array}{ll}
x(n+1)= \cfrac{\lambda x(n)}{1+x(n)}\displaystyle e^{-y(n)(1+\alpha y(n))}\\[1ex]
y(n+1)=\beta x(n)\left(1-e^{-y(n)(1+\alpha y(n))}\right).
\end{array}\right.
\end{equation}

We first observe that solutions of \eqref{rmodel} remain
nonnegative and are bounded for $n\geq 0$. The trivial steady
state $E_0=(0,0)$ exists for all feasible parameters and the
Jacobian matrix of \eqref{gen} at $(x,y)$ is given by
\begin{equation}\label{J}J(x,y)=\left(\begin{array} {cc}
\cfrac{\lambda e^{-y(1+\alpha y)}}{(1+x)^2} & -\cfrac{\lambda xe^{-y(1+\alpha y)}(1+2\alpha y)}{(1+x)}\\[2ex]
\beta \bigg(1-e^{-y(1+2\alpha y)}\bigg) & \beta xe^{-y(1+\alpha
y)}(1+2\alpha y)
\end{array}
\right).
\end{equation}
At $E_0=(0,0)$, $J(E_0)=\left(\begin{array}
{cc}
\lambda  & 0\\[1ex]
0 & 0
\end{array}
\right)$ and hence $E_0$ is asymptotically stable if $\lambda<1$
and it is a saddle point with the stable manifold lying on the
nonnegative $y$-axis if $\lambda>1$.

\medskip

Notice  $x(n+1)\leq \lambda x(n)/(1+ x(n))$ for $n\geq 0$ implies
$\lim_{n\rightarrow\infty}x(n)=0$ if $\lambda \leq 1$. Hence
$\lim_{n\rightarrow\infty}y(n)=0$ and $E_0$ is globally
attracting, and we have the following result.

\medskip

\noindent{\bf Proposition 2.1} {\em If $\lambda \leq 1$, then
$E_0=(0,0)$ is globally asymptotically stable for \eqref{rmodel}.}

\medskip

Proposition 2.1 implies that if the intrinsic growth rate
$\lambda$ of the prey population  is not greater than one, then
the prey population goes extinct and so does the predator
population.

We assume $\lambda>1$ for the remainder discussion so that the
prey population can persist in the absence of predator. It follows
that \eqref{gen} has another boundary steady state
$$
E_1=(\bar x, 0), \mbox{ where } \bar x=\lambda-1>0.
$$
Notice that $\bar x$ can be viewed as the carrying capacity of the
prey population. The Jacobian matrix of \eqref{rmodel} evaluated
at $E_1$ is
$$
J(E_1)=\left(\begin{array} {cc} 1+\lambda
\bar xg_0^\prime(\bar x) & *\\[1ex]
0 & \beta  \bar x
\end{array}
\right), \mbox{ where * is an unimportant term.}
$$
 Observe that
$0<1+\lambda \bar xg_0^\prime(\bar x)<1$. Therefore $E_1$ is
asymptotically stable if $\beta \bar x<1$ and it is a saddle point
with its stable manifold lying on the positive $x$-axis if $\beta
\bar x>1$.

Since there are only two boundary steady states for which their
stability is known, we prove that system \eqref{rmodel} is
uniformly persistent when $E_1$ is unstable. That is, there exists
$\eta>0$ such that $\liminf_{n\rightarrow\infty}x(n)\geq\eta$ and
$\liminf_{n\rightarrow\infty}y(n)\geq\eta$ for all solutions of
\eqref{rmodel} with $x(0)>0$ and $y(0)>0$.   Our proof is based on
the boundary dynamics of \eqref{rmodel} using Theorem 4.1 of
\cite{so}.

\medskip

\noindent{\bf Theorem 2.2} {\em Let $\lambda>1$ and $\beta  \bar
x>1$. Then system \eqref{rmodel} is uniformly persistent.}

\begin{proof} Let $Y$ be the boundary of the nonnegative
coordinate plane $\mathbb{R}_+^2$. Then $\mathbb{R}_+^2\backslash
Y$ is positively invariant for system \eqref{rmodel}. Clearly
solutions $(x(n),y(n))$ of \eqref{rmodel} satisfy
$\limsup_{n\rightarrow\infty}x(n)\leq\bar x$ and
$\limsup_{n\rightarrow\infty}y(n)\leq\beta\bar x$, and hence
system \eqref{rmodel} has a global attractor $X$. The only
invariant sets in $Y$ are $\{E_0\}$ and $\{E_1\}$, where $E_i\in
X$ for $i=0,1$. Applying Theorem 4.1 of \cite{so}, we need to
verify that each $\{E_i\}$ is isolated in $X$ and the stable set
of $E_i$ is contained in $Y$. Since $X$ is closed in
$\mathbb{R}_+^2$, it is sufficient to show that $\{E_i\}$ is
isolated in $\mathbb{R}_+^2$ for $i=1, 2$.

If $\{E_0\}$ is not isolated in $\mathbb{R}_+^2$, then for any
$\epsilon>0$ there exists a compact invariant set $M_0$ in
$\overline{B(E_0,\epsilon)}\bigcap \mathbb{R}_+^2$ such that
$\{E_0\}$ is a proper subset of $M_0$, where $B(E_0, \epsilon)$
denotes the $\epsilon$-ball centered at $E_0$. Since $\lambda>1$,
we can choose $\epsilon>0$ so that
$\cfrac{\lambda}{1+\epsilon}e^{-\epsilon(1+\alpha \epsilon)}>1$.
Let $x^*=\sup\{x:(x,y)\in M_0\}$. Then $0<x^*\leq\epsilon$ and
there exists $y^*\leq\epsilon$ such that $(x^*,y^*)\in M_0$. Let
$x(0)=x^*$ and $y(0)=y^*$. Then $(x(n),y(n))\in M_0$ for $n\geq 0$
and
$$x(1)\geq
x(0)\cfrac{\lambda}{1+\epsilon}e^{-\epsilon(1+\alpha
\epsilon)}>x(0)=x^*.$$ We obtain a contradiction and conclude that
$\{E_0\}$ is isolated in $\mathbb{R}_+^2$.

Suppose now $\{E_1\}$ is not isolated in $\mathbb{R}_+^2$. Then
for any $\epsilon>0$ there exists a compact invariant set $M_1$ in
$\overline{B(E_1,\epsilon)}\bigcap\mathbb{R}_+^2$ with $\{E_1\}
\varsubsetneq M_1$. We choose $\epsilon>0$ such that $\beta(\bar
x-\epsilon)(1-\cfrac{\epsilon}{2})>1$ and $a\epsilon<2$. Then
$\beta(\bar x-\epsilon)>1$. Let $(x(0),y(0))\in M_1$ with
$y(0)>0$. Then $y(n)>0$  and $(x(n),y(n))\in M_1$ for $n\geq 0$.
It follows that $y(n+1)=\beta x(n)(1-e^{-y(n)(1+\alpha
y(n))})>\beta(\bar x-\epsilon)(1-e^{-y(n)})$ for $n\geq 0$.
Consider $z(n+1)=\beta(\bar x-\epsilon)(1-e^{-z(n)})$ with
$z(0)=y(0)$. Since $\beta(\bar x-\epsilon)>1$, the scalar equation
has a unique interior steady state $\bar z$ such that
$\lim_{n\rightarrow\infty}z(n)=\bar z$ if $z(0)>0$. It follows
that $\liminf_{n\rightarrow\infty}y(n)\geq \bar z$. We claim $\bar
z>\epsilon$. Indeed, $\beta(\bar
x-\epsilon)(1-e^{-\epsilon})/\epsilon=\beta(\bar
x-\epsilon)\bigg(1-\cfrac{\epsilon}{2}+\cfrac{\epsilon^2}{6}-\cfrac{\epsilon^3}{24}+....\bigg)>\beta(\bar
x-\epsilon)(1-\cfrac{\epsilon}{2})$ since $4>\epsilon$. Therefore,
$\beta(\bar x-\epsilon)(1-e^{-\epsilon})/\epsilon>1$ and
$\epsilon<\bar z$ is shown. We then have
$\liminf_{n\rightarrow}y(n)\geq\bar z>\epsilon$ and obtain a
contradiction. Consequently, $\{E_1\}$ is isolated in
$\mathbb{R}_+^2$.  It is straightforward to see that the stable
sets of $E_0$ and $E_1$ lie on $Y$ and hence system \eqref{rmodel}
is uniformly persistent by Theorem 4.1 of \cite{so}. \end{proof}

Theorem 2.2 indicates that both prey and predator populations can
coexist if $\beta \bar x>1$, where $\beta  \bar x$ can be viewed
as the maximal reproductive number of the predator since the prey
population is stabilized at the carrying capacity $\bar x$ level.
The parameter $\alpha$ plays no role in the sufficient condition
for coexistence. On the other hand if $\beta \bar x<1$, then since
$E_1=(\bar x, 0)$ is asymptotically stable, the system is not
uniformly persistent.

To study the effects of cooperative hunting, we need to understand
the dynamics of the model when there is no cooperative hunting in
the predator. The dynamics of such a  model are given in Theorem
3.1 of \cite{jang2006} and are restated as follows.

\medskip

\noindent{\bf Proposition 2.3} {\em Let $\alpha=0$. The trivial
steady state $E_0=(0,0)$ is globally asymptotically stable if
$\lambda\leq 1$. If $\lambda>1$, then \eqref{rmodel} has another
boundary state $E_1=(\bar x,0)$ which is globally asymptotically
stable if $\beta \bar x <1$.  If $\beta \bar x >1$, then
\eqref{rmodel} has a unique interior steady state and the system
is uniformly persistent.}

\medskip

 Although
it is not proved analytically in \cite{jang2006}, it is observed
that when $\alpha=0$ the unique interior steady state loses its
stability via a Neimark-Sacker bifurcation as $\beta$ increases.

\section{Interior steady states for $\alpha>0$}
\setcounter{equation}{0}

In this section, we study existence and  number of  interior
steady states. These are achieved by investigating geometry of the
isoclines.

The nontrivial $y$-isocline is given by
\begin{equation}\label{h}
x=h(y):=\cfrac{y}{\beta(1-e^{-y(1+\alpha y)})} >0 \mbox{ with }
h(0)=\cfrac{1}{\beta } \mbox{ and }  \ h(\infty)=\infty.
\end{equation}
For simplicity, we introduce a new notation
\begin{equation}\label{sym}
\diamondsuit=y(1+\alpha y).
\end{equation}
Then
\begin{equation}\label{g}
h^\prime(y)=\cfrac{e^{-\diamondsuit}
g(y)}{\beta(1-e^{-\diamondsuit })^2}  \mbox{ with }
h^\prime(0)=\cfrac{1-2\alpha}{2 \beta },
\end{equation}
where $g(y)=e^{\diamondsuit}-1-y(1+2\alpha y)$. It is easy to
check that
$$
 g^{\prime}(y)= \{e^{\diamondsuit}(1+2\alpha y)-(1+4\alpha y) \} \mbox{ and }
  g^{\prime\prime}(y)=\{ e^{\diamondsuit}(2\alpha+ (1+2\alpha y)^2)-4\alpha \}.
 $$
 Then $ g(0)=0$, $ g^\prime(0)=0$,
 $  g^{\prime\prime}(0)=(1-2\alpha)$ and  $ g^{\prime\prime\prime}(y)>0$ for $y\geq 0$.
 If $2\alpha\leq 1$, we have from $ g(0)=0, g^\prime(0)=0$ and  $  g^{\prime\prime}(0)=(1-2\alpha) \geq 0 $
 that  all $g^{\prime\prime}(y), g^{\prime}(y) $ and $g(y)>0$ for $y>0$.  So $h^\prime(y)>0$ for $y>0$ as sign $(h^\prime) =$ sign $(  g)$ by \eqref{g}. Thus
\begin{equation}\label{monotone h}
h \uparrow \mbox { strictly from } h(0)=\cfrac{1}{\beta }  \mbox{
to } \ h(\infty)=\infty   \mbox{  if }  2\alpha\leq 1.
\end{equation}
If $2\alpha>1$ then $ g^{\prime\prime}(0)<0$ and
$h^{\prime}(0)<0$. Using the same argument we can show that there
exists a unique critical point $\bar y>0$ such that
$h^\prime(y)<0$ on $[0, \bar y)$ and $h^\prime(y)>0$ on $(\bar y,
\infty)$. Hence, while $h>0$ on $[0, \infty)$,
\begin{equation}\label{nonmonotone h}
h \downarrow \mbox {strictly  on } [0,  \bar y]  \mbox{ and then }
h \uparrow  \mbox { strictly on }  [\bar y, \infty) \mbox{ if } <
2\alpha>1.
\end{equation}

The non-trivial $x$-isocline is given by
\begin{equation}\label{f}
x=f(y):=\cfrac{\lambda e^{-y(1+\alpha y)}-1}{k}
\end{equation}
with $ f^\prime(y)=-\lambda e^{\diamondsuit}(1+2\alpha y)<0$.
Define
\begin{equation}\label{y_c}
y_c=\cfrac{-1+\sqrt{1+4\alpha\ln\lambda}}{2\alpha}
\end{equation}
by solving $f(y)=0$, i.e., $e^{\diamondsuit}=\lambda$. Then
\begin{equation}\label{monotone f}
  f \downarrow \mbox { strictly on } [0,  \infty) \mbox{ with } f(0)=\bar x=\lambda-1>0   \mbox{ and } f(y_c)=0.
\end{equation}
For the existence of interior steady states, we are only concerned
with $y \in (0,y_c)$ since the $x$ component of the steady state
would be negative if $y>y_c$.


For $2\alpha\leq 1$, it follows from \eqref{monotone h} and
\eqref{monotone f} that system \eqref{rmodel} has either zero or
one interior steady state depending on whether $1/\beta   = h(0)
\geq f(0)=\bar x$, i.e., whether  $\beta  \bar x \leq 1$.  This
proves part of the following result.

\medskip

\noindent{\bf Theorem 3.1} {\em Let $\lambda>1$. Then system
\eqref{rmodel} has a unique interior steady state for $\beta  \bar
x>1$. If $0<2\alpha\leq \cfrac{3\lambda-1}{\lambda-1}$ and $\beta
 \bar x\leq 1$ then \eqref{rmodel} has no interior steady state. If
$2\alpha> \cfrac{3\lambda-1}{\lambda-1}$, then \eqref{rmodel} has
exactly one interior steady state in case $\beta \bar x \geq 1$
and at most two interior  steady states in case $\beta  \bar x
<1$.}

\begin{proof} It suffices to consider $2\alpha>1$. By \eqref{nonmonotone h}, the
$y$-isocline is no longer monotone and therefore the analysis is
different  from that of the case $2\alpha\leq 1$. Instead of
analyzing the convexity property of both isoclines, we adopt a
different approach. Setting the two nontrivial isoclines equal,
$h(y)=f(y)$, it leads to solve the following system on $y \in (0,
y_c)$
\begin{equation}\label{intersection}\left\{\begin{array}{ll}
z=w(y):=\beta\bigg(\lambda e^{-\diamondsuit}-1\bigg)\bigg(1-e^{-\diamondsuit}\bigg) \\[1ex]
z=q(y):=y.
\end{array}\right.
\end{equation}
 In fact, solving \eqref{intersection} is equivalent to solve for interior steady states of \eqref{rmodel}.
  We need to study the convexity property of the function $w$ defined in the first equation. Notice
\begin{equation}\label{w}
w(0)=0=w(y_c), \ w(y)>0 \mbox{ on }  (0, y_c)
\end{equation}
and
\begin{equation}\label{wprime}
w^\prime(y)=\beta\bigg(-(\lambda+1)e^{-\diamondsuit}(1+2\alpha
y)+2\lambda e^{-2\diamondsuit}(1+2\alpha y)\bigg).
\end{equation}
Observe that $w^\prime(y)=0$ has a unique positive solution $y_1 <y_c$ satisfying $e^{ \diamondsuit}=\cfrac{2\lambda }{\lambda+1} <\lambda$. Hence, $w(y_1)= \max w $ and
\begin{equation}\label{y1}
w \uparrow  \mbox{ strictly on }[0, y_1) \mbox{ and } w \downarrow   \mbox{ strictly on } (y_1, \infty).
\end{equation}
Define $u(y)=-(\lambda+1)e^\diamondsuit[2\alpha-(1+2\alpha
y)^2]+2\lambda[2\alpha-2(1+2\alpha y)^2]$. Then
 \begin{equation}\label{wprimeprime}
w^{\prime\prime}(y)=\beta e^{-2\diamondsuit} u(y).
\end{equation}
Using  $u(0)=-(\lambda+1)(2\alpha-1)+4\lambda(\alpha-1) =2\alpha
(\lambda-1)-(3\lambda-1)$, we have
\begin{equation}\label{u0}
u(0)>  0 \mbox{ if and only if } 2\alpha >
\cfrac{3\lambda-1}{\lambda-1}.
\end{equation}
 By a direct computation,
 \begin{equation}\label{uprime}
 u^\prime(y)=(1+2\alpha y)v(y),
 \end{equation}
 where
$v(y)=(\lambda+1)e^\diamondsuit[(1+2\alpha
y)^2+2\alpha]-16\lambda\alpha .$ Apparently  $v^\prime(y)>0$ and
$v(\infty)= \infty$.  Using $\lambda >1$ and $1<2\alpha$, we can
show that
$$
u^\prime(0)=v(0)=  (\lambda+1)(1 +2\alpha)-16\lambda\alpha  \leq 2
\lambda (1-6 \alpha )< 0.
$$
Hence, there exists $y_2 >0$ such that $v<0$ on $[0, y_2)$ and
$v>0$ on $(y_2, \infty)$. It follows from \eqref{uprime} that
\begin{equation}\label{u}
 u \downarrow  \mbox{ on }[0, y_2) \mbox{ and } u \uparrow   \mbox{ on } (y_2, \infty).
 \mbox{ In particular, }  u(y_2)= \min u.
 \end{equation}
 Using $u(\infty)= \infty$, \eqref{wprimeprime} and \eqref{u0},   there exists $y_0>0$ such that
  \begin{equation}\label{wprime2a}
 w^{\prime\prime} <0 \mbox{ on } (0, y_0) \mbox{ and } w^{\prime\prime} >0    \mbox{ on } (y_0, \infty) \mbox{ for }
 1< 2\alpha \leq \cfrac{3\lambda-1}{\lambda-1}.
 \end{equation}
Note that $y_1 \leq y_0$ as $ w^{\prime\prime}(y_1) \leq 0$ by \eqref{y1}.

In case $2\alpha >\cfrac{3\lambda-1}{\lambda-1}$, we have $u(0)> 0
$ by \eqref{u0}. We claim $u(y_2)<0$. Otherwise, $
w^{\prime\prime} \geq 0$ by \eqref{u} and \eqref{wprimeprime},
which leads to a contradiction to \eqref{w} and \eqref{y1}.
 Using $ w^{\prime\prime}(y_1) \leq 0$ and $u(\infty)= \infty$ again, there exist $0< y_3 < y_4$ with $y_1 \in [y_3, y_4]$ such that
\begin{equation}\label{wprime2b}
 w^{\prime\prime} >0 \mbox{ on } [0, y_3) \cup (y_4, \infty) \mbox{ and } w^{\prime\prime} <0    \mbox{ on } (y_3, y_4)
 \mbox{ for }  2\alpha > \cfrac{3\lambda-1}{\lambda-1}.
\end{equation}

Now we are ready to prove the conclusions of the theorem. Consider
first the case $1< 2\alpha \leq \cfrac{3\lambda-1}{\lambda-1}$.
Note that by \eqref{wprime},
\begin{equation}\label{wprime0}
w^\prime(0)=\beta(\lambda-1)>0.
\end{equation}
If $y_c \leq y_0$,  then $w$ is concave on the whole interval $[0,y_c]$ by \eqref{wprime2a}. Since $w(0)=q(0)=0$,
it follows easily that system \eqref{wprime} has none or one intersection on $(0,y_c)$ depending on whether
\begin{equation}\label{k}
 1=q^\prime(0)\geq w^\prime(0)=\beta(\lambda-1),  \mbox{ i.e., }  \beta \bar x \leq 1.
\end{equation}

In case $y_c > y_0$,  $w$ is concave on the interval $[0,y_0]$ and convex on $[y_0,y_c]$. We will count the number of intersection points on each interval above.
Under \eqref{k}, there is no interior intersection on
 $(0,y_0]$ as before. There is no  intersection on
 $(y_0,y_c)$ either due to $q(y_0) > w(y_0)$ or by  \eqref{y1},
 \begin{equation}\label{qw}
 q \uparrow  \mbox{ strictly on } [y_0,y_c]  \mbox{ and } w \downarrow \mbox{ strictly  on } [y_0,y_c].
 \end{equation}
 If $ \beta \bar x > 1$, then \eqref{wprime0} implies that $k= q^\prime(0)< w^\prime(0)$. We need to compare $ q(y_0)$ with $ w(y_0)$. By \eqref{wprime2a} again,
 $$
 q(y_0) \geq w(y_0) \mbox{ if and only if  system  \eqref{intersection} has a solution on } (0,y_0].
 $$
On the other hand,   \eqref{qw} implies that
 $$
 q(y_0) < w(y_0) \mbox{ if and only if  system }  \eqref{intersection} \mbox{ has a solution on } (y_0,y_c].
 $$
 Summing up, system  \eqref{intersection} has exactly one solution on $(0,y_c)$ for $ \beta \bar x > 1$.

 It remains to consider the case $ 2\alpha > \cfrac{3\lambda-1}{\lambda-1}$. It suffices  to show
 \begin{equation}\label{bax1}
 \mbox{ if }  \beta \bar x \geq 1,  \mbox{ then } (3.9) \mbox{ has exactly one solution on } (0,y_c)
 \end{equation}
 and
  \begin{equation}\label{bax2}
 \mbox{ if }  \beta \bar x < 1,  \mbox{ then } (3.9) \mbox { has at most two solution on } (0,y_c).
 \end{equation}
 In view of \eqref{wprime2b}, we need to compare $y_4$ with $y_c$. In case  $y_4 < y_c$, we will count separately the number of intersection points
 on $(0,y_3]$, $(y_3,y_4]$ and $(y_4,y_c)$. Denote them by $N_1, N_2$ and $N_3$ respectively. First note that by \eqref{wprime2b} and \eqref{y1},
 \begin{equation}\label{N2}
\mbox{  if } q(y_3) <w(y_3)     \mbox{ then } N_2+N_3 =1.
\end{equation}
In fact  $N_2+N_3 =1+0 $ or $0+1$ depending on whether  $q(y_4)  \geq w(y_4)$.

If $ \beta \bar x \geq 1$, then \eqref{wprime0} implies
$q^\prime(0) \leq w^\prime(0) $. Further,  $q(y_3) <w(y_3)$  and
$N_1=0$ by \eqref{wprime2b}. Hence \eqref{bax1} follows from
\eqref{N2}.

 If $ \beta \bar x < 1$, then $q^\prime(0) > w^\prime(0) $ and \eqref{wprime2b} implies $N_1=0$ or $1$ depending on whether $q(y_3)>w(y_3)$.
 In case $q(y_3) >w(y_3)$ we may repeat the same arguments as for the case $ 2\alpha \leq \cfrac{3\lambda-1}{\lambda-1}$ discussed above.
 Depending on whether $q(y_4) \geq w(y_4)$,  $N_2 +N_3$ can be $0+0$, $1+0$, $2+0$ or $1+1$. The last case happens when $q(y_4) < w(y_4)$.  So $N_1+N_2+N_3 \leq 2$
 as claimed by \eqref{bax2}, which still holds if  $q(y_3) =w(y_3)$, except that  $N_1$ is increased by 1 and $N_2+N_3 $ decreased by 1.
 If $q(y_3) <w(y_3)$, we have $N_1= 1$ and $ N_2+N_3 =1$ by \eqref{N2}. Therefore, \eqref{bax2} follows.

Note that in the counting of $N_3$ we only use the fact that on $(y_4,y_c)$,
$$
q \uparrow \mbox{ strictly and } w \downarrow \mbox{ strictly},
$$
which holds for $(y_1,y_c)$ as well by \eqref{y1}. Therefore the
same argument also works for the remaining case of $y_4 \geq y_c$
by counting the number of intersection points on  $(0,y_3]$,
$(y_3,y_1]$ and $(y_1,y_c)$ separately. The detail is omitted.
\end{proof}

\medskip

Recall that $\beta \bar x$ is the maximal reproductive number of
the predator since the prey population is stabilized at its
carrying capacity $\bar x$. If this reproductive number exceeds
one, then the predator-prey interaction can support a unique
coexisting steady state. If this reproductive number is smaller
than one and the degree of cooperation is also small, then Theorem
3.1 implies that the predator-prey interaction has no coexisting
steady state while the interaction may support two coexisting
steady states if the degree of predator cooperation is large.

When $\lambda>1$, $2\alpha> \cfrac{3\lambda-1}{\lambda-1}$ and
$\beta < \cfrac{1}{\bar x}$, the number of interior steady states
is not clear from Theorem 3.1. To clarify this issue, we let the
$y$-isocline \eqref{h}  vary with $\beta$  while fixing the
$x$-isocline \eqref{f}. Note that by  \eqref{y_c} and
\eqref{monotone f}, we are only concerned with the interval
$(0,y_c)$.    Rewrite the  $y$-isocline as
$x=h_{\beta}(y)=\cfrac{y}{\beta(1-e^{-y(1+\alpha y)})} $ to
emphasize its dependence on $\beta$. For fixed $y>0, h_{\beta}(y)$
decreases strictly to $0$ as $\beta$ increases to $\infty$.
Therefore on the interval $ [0,y_c)$,
\begin{equation}\label{bfamily}
  \{  x=h_{\beta}(y) :   \beta >0 \} \mbox{ forms a family of nonintersecting curves}.
\end{equation}
We let $\beta $ increase from $0$. Apparently, both isoclines do not intersect for $\beta $ small. At a certain $\beta_* $,
both isoclines become tangent to each other. The tangent point is unique.
Otherwise by increasing $\beta $ a little over $\beta_* $, two isoclines will have four intersection points, which is contrary to Theorem 3.1.
 Denote the tangent point by  $E_*=(x_*, y_*)$. In particular, \eqref{T} holds. For $\beta_* <\beta < 1/\bar x$,
 we have $ 1/\beta =h_{\beta}(0) > f(0) = \bar x, h_{\beta}(y_*) <h_{\beta_*}(y_*)= f(y_*)$ and  $h_{\beta}(y_c) > 0=f(y_c)$.
 Two isoclines will have exactly two interior steady states in view of Theorem 3.1. Denote them by $E_i^*=(x_i^*,y_i^*)$, $i=1,2$,
 with $0<y_1^*<y_*<y_2^*<y_c$. Note that as $\beta$ increases, $E_1^* $ moves to the left and $E_2^* $ moves to the right along the
 the $x$-isocline. When $\beta \uparrow 1/\bar x$, $E_1^* $ becomes $(\bar x,0)$, no longer an interior steady state,
 and $ \lim E_2^* $  exists. Denote it by $E^*=(x_e, y_e).$  This is consistent with  Theorem 3.1 which indicates that there is exactly
  one interior steady state for $\beta  \geq 1/\bar x$.  The
discussion is summarized as follows.

\medskip

 \noindent{\bf Theorem 3.2} {\em Let $\lambda>1$ and
$2\alpha>\cfrac{3\lambda-1}{\lambda-1}$. Then there exists a
unique $\beta_*>0$ such that \eqref{rmodel} has no interior steady
state if $\beta<\beta_*$.  When $\beta=\beta_*$, \eqref{rmodel}
has a unique interior steady state $E_*=(x_*, y_*)$ at which both
isoclines intersect tangentially and $y_*\in (0, y_c)$ is uniquely
determined by
\begin{equation}\label{T}\left\{\begin{array}{ll}
f(y)= h(y)\\[1ex]
f^\prime(y)=h^\prime(y).
\end{array}\right.
\end{equation}
For $\beta\in(\beta_*, \cfrac{1}{\bar x})$, there are two interior
steady states $E_i^*=(x_i^*,y_i^*)$, $i=1,2$, with $y_1^*<y_2^*$.
Moreover, there is a unique interior steady state
 if $\beta \geq \cfrac{1}{\bar x}$. Denote it by $E^*=(x_e,y_e)$ for $\beta =\cfrac{1}{\bar x}$ and denote it by $E^*=(x^*,y^*)$ for   $\beta > \cfrac{1}{\bar x}$.
 Hence, we have
 \begin{equation}\label{order}
 0< y_1^*< y_* <y_2^*< y_e < y^*< y_c.
 \end{equation}   }


Theorem 3.2 provides a criterion in terms of $\beta$ for which the
predator-prey interaction can support two coexisting steady
states, where $\beta$ is the predator conversion for each prey
consumed. If the maximal reproductive number of the predator is
smaller than one, then system \eqref{rmodel} has two coexisting
steady states if $\beta$ is larger than the critical value
$\beta_*$.

 Note that \eqref{bfamily} implies that for $\beta \in
( \cfrac{1}{\bar x}, \infty)$,
$$
x^* \downarrow \mbox{ strictly and } y^* \uparrow \mbox{ strictly as } \beta \mbox{ increases.}
$$
This property holds as long as $\lambda >1$ and  $\beta >
\cfrac{1}{\bar x}$ for which Theorem 3.1 guarantees that there is
exactly one interior steady state $E^*=(x^*,y^*)$.

\bigskip

Figure 1 plots the two isoclines under different scenarios. In (a)
$\lambda=5$ and $\alpha=1/2.1$ so that $1\geq 2\alpha$. Two
$\beta$ values of $0.21$ and $0.525$ are chosen to illustrate the
nonexistence and existence of a unique interior steady state
respectively. In (b) $\lambda=15$, $\alpha=6/5$ so that
$1<2\alpha$ and $2\alpha-\cfrac{3\lambda-1}{\lambda-1}=-0.743<0$.
Two $\beta$ values of $0.05$ and $0.125$ are chosen to show
respectively the nonexistence and existence of a unique interior
steady state. In (c) $\lambda=10$ and $\alpha=15$, where
$2\alpha-\cfrac{3\lambda-1}{\lambda-1}\approx 26.78>0$ and
$2\alpha>1$. The two curves are tangent to each other at
$\beta=\beta_*\approx 0.066502$, and there is a unique positive
intersection when $\beta=\cfrac{1}{\bar x}\approx 0.11$. It is
clear that the system has two interior steady states if $\beta$ is
in $(0.066502, 0.11)$, there is a unique interior steady state if
$\beta>0.11$  and there is no interior steady state if
$\beta<\beta_*=0.066502$.

\section{Stability and dynamics of the model}
\setcounter{equation}{0} To study  asymptotic dynamics of
\eqref{rmodel} for $\lambda>1$ and $\alpha>0$, we  separate our
discussion into two cases: $2\alpha\leq 1$, and $2\alpha>1$. Since
local stability can be determined by the Jury conditions, we first
provide a result on the determinant of the Jacobian matrix $J$ at
an interior steady state $E=(x,y)$. Using \eqref{h} and \eqref{f}
we may rewrite $J$ as follow

\begin{equation}\label{J1st}
J(E) =\left(\begin{array} {cc} \cfrac{1}{\lambda e^{-\diamondsuit}} & -x(1+2\alpha y) \\[2ex]
 \cfrac{y}{x} &   \cfrac{y(1+2\alpha y)e^{-\diamondsuit}}{1-e^{-\diamondsuit}}
\end{array}
\right) \mbox{ with } x= f(y).
\end{equation}

Recall $\diamondsuit =y(1+ \alpha y)$ by \eqref{sym}. Following
Theorem 3.2, we will let $E$ vary along  the $x$-isocline and thus
$J $ in \eqref{J1st} becomes a function of $y$.

\medskip

\noindent{\bf Lemma 4.1} {\em Let $\lambda>1$. With $y_c$ defined in \eqref{y_c}, $det(J)$ is a
strictly increasing function of $y$ with  $det(J)|_{y=0} <1$ and $det(J)|_{y=y_c} >1$.
Hence, there exists a unique $y_d\in (0, y_c)$ such that  $det(J)|_{y=y_d} =1$.}

\begin{proof} By \eqref{J1st},
\begin{equation}\label{det}
det(J)= \cfrac{y(1+2\alpha y)}{\lambda[1-e^{-\diamondsuit}]}+
y(1+2\alpha y).
\end{equation}
Since $\{y(1+2\alpha y)\}^{\prime} =(1+4\alpha y)>0$ for $y>0$,
$$
\cfrac{d }{d y}det(J) \geq  \cfrac{ae^{-\diamondsuit} \{
(e^\diamondsuit-1)(1+4\alpha y)- y(1+2\alpha y)^2\}  }{\lambda
[1-e^{-\diamondsuit}]^2}  >0
$$
as the term inside the bracket  $ \geq  \diamondsuit(1+4\alpha y)-
y(1+2\alpha y)^2=\alpha y^2 >0$.

By \eqref{det}, $det(J)|_{y=0^+}=\cfrac{1}{\lambda}<1$. Since $y_c
$ satisfies $e^\diamondsuit=\lambda$, $y_c(1+2\alpha y_c)
>y_c(1+\alpha y_c) = \ln\lambda$ and thus $det(J)|_{y=y_c} \geq
\cfrac{\lambda \ln\lambda}{\lambda-1}>1$ as it is an increasing
function   with $\lim_{ \lambda \downarrow 1
}\cfrac{\lambda\ln\lambda}{\lambda-1}=1$. The conclusion
follows.\end{proof}

\medskip

Figure 2(a) plots $det(J)$ as a function of $y$ with $\lambda=10$,
$\beta=6.3/20$ and $\alpha=5/2.1$, where $det(J(y))$ crosses the
horizontal line $det(J)=1$ at $y=y_d$.

\subsection{Dynamics of the model when $2\alpha\leq 1$}

Under this assumption, Theorem 3.1 implies that \eqref{rmodel} has
a unique interior steady state $E=(x, y)$ if $\beta  \bar x>1$,
and there is no interior steady state if $\beta \bar x<1$. Our
analysis is more complete in this parameter regime. In particular,
$E_1=(\bar x,0)$ is globally asymptotically stable if $\beta \bar
x<1$, which is similar to the system with $\alpha=0$ as
illustrated in Proposition 2.3.

\medskip

\noindent{\bf Theorem 4.2} {\em Let $\lambda>1$, $ 2\alpha\leq 1$
and $\beta \bar x<1$. Then $E_1=(\bar x, 0)$ is globally
asymptotically stable in $\{(x,y)\in\mathbb{R}_+^2: x>0\}$.}

\begin{proof} Under the given assumptions, $E_1$ is locally asymptotically
stable and there is no interior steady state. It is
straightforward to verify that
\begin{equation}\label{ineq}
1-e^{-y(1+\alpha y)}<y \mbox{ for } y>0.
\end{equation}
Indeed, letting $\tilde h(y)=1-e^{-y(1+\alpha y)}-y$, we obtain
$\tilde h(0)=0$, $\tilde h^\prime(y)=e^{-y(1+\alpha y)}(1+2\alpha
y)-1$, $\tilde h^\prime(0)=0$  and $\tilde
h^{\prime\prime}(y)=e^{-y(1+\alpha y)}\big(2\alpha-(1+2\alpha
y)^2\big)<0$ since $1\geq 2\alpha$. Hence, $\tilde h(y) <0 $ for
$y>0$ and \eqref{ineq} is verified.

By \eqref{rmodel}, $x(n+1)\leq\cfrac{\lambda x(n)}{1+x(n)}$. This
implies $\limsup_{n\rightarrow\infty}x(n)\leq \bar x$. Thus for
any $\epsilon>0$, $x(n)<\bar x+\epsilon$ for $n $ large. We choose
$\epsilon>0$ such that $\beta (\bar x+\epsilon)<1$. Let
$(x(0),y(0))\in \mathbb{R}_+^2$ with $x(0)>0, y(0) \geq 0$. By
\eqref{rmodel} and \eqref{ineq}, $y(n+1)<\beta x(n) y(n)<\beta
(\bar x+\epsilon)y(n)$ for $n $ large.  Hence,
$\lim_{n\rightarrow\infty}y(n)=0$.
 Then for any $\eta>0$, $x(n+1)\geq\cfrac{(\lambda -\eta) x(n)}{1+x(n)}$   for $n $ large by  \eqref{rmodel}.
 Consequently, $\liminf_{n\rightarrow\infty}x(n)\geq \bar x$ and thus $\lim_{n\rightarrow\infty}x(n)=\bar x$. Therefore,  $E_1$ is
globally attracting in $\{(x,y)\in\mathbb{R}_+^2: x>0\}$ as claimed. \end{proof}

\medskip

Theorem 4.2 states that if the intrinsic growth rate $\lambda$ of
the prey population is greater than one, the degree of predator's
cooperative hunting is small, $2\alpha\leq 1$, and the maximal
reproductive number of the predator  is less than one, $\beta \bar
x<1$, then the predators will go extinct and the prey population
will stabilize at its carrying capacity $\bar x$.

Let $\beta \bar x>1$. By \eqref{monotone h} and Theorem 3.1,
\eqref{rmodel} has a unique interior steady state $E=(x, y)$ with
$\cfrac{1}{\beta }<x<\bar x$ and $0<y<y_c$. Moreover,
\eqref{bfamily} implies that
\begin{equation}\label{b_monotone}
x \downarrow \mbox{ strictly  and } y \uparrow \mbox{ strictly as
} \beta \mbox{ increases  from }  \cfrac{1}{\bar x} \mbox{ to }
\infty.
\end{equation}
In fact, $(x,y)$ converges to $(0,y_c)$ along the x-isocline as $ \beta \to   \infty$.
Rewrite the Jacobian matrix \eqref{J1st} at $E$ as
\begin{equation}\label{J2nd}
J(E )=\left(\begin{array} {cc} \tilde a & -\tilde b \\[1ex]
\tilde c & \tilde d
\end{array}
\right),
\end{equation}
where $\tilde a>0, \ \tilde b>0, \ \tilde
c>0$ and $\tilde d>0$. Clearly
\begin{equation}\label{Jury}
-tr(J)< 1+det(J)
\end{equation}
as $det(J)=\tilde a \tilde d +\tilde b \tilde c  >0$ and $tr(J)= \tilde
a  + \tilde d >0$. By the Jury
conditions \cite{allen}, the local stability of $E$ is determined by
 $det(J)$, which increases strictly by Lemma 4.1, and
 \begin{equation}\label{V}
V(y):=1+det(J)-tr(J),
\end{equation}
which behaves well at the present case $1\geq 2\alpha$. Indeed,
using \eqref{J2nd}, we have
\begin{equation}\label{V1}
V(y)=(1-\tilde a)(1-\tilde d)+\tilde b\tilde c > \tilde b\tilde c >0
\end{equation}
as $\tilde a=1/(1+x)<1$ by \eqref{f} and $\tilde
d=\cfrac{y(1+2\alpha y)}{e^{\diamondsuit}-1} < 1$ due to
$e^{\diamondsuit}-1\geq y(1+\alpha y)+\cfrac{y^2(1+\alpha
y)^2}{2}> y(1+\alpha y+\cfrac{y}{2})\geq y(1+2\alpha y)$ for
$y>0$. In the last inequality the assumption $1\geq 2\alpha$ is
used.



\medskip

Since  magnitude of the interior steady state is monotone with
respect to $\beta$ by \eqref{b_monotone} and $\beta$ is the prey
conversion to predator, we use $\beta$ as the bifurcation
parameter. Let $\beta_d>0$ be the corresponding $\beta$ value of
$y_d$ given in Lemma 4.1. The following result follows from Lemma
4.1 and \eqref{V1}.   Let $C^*$ be defined in \eqref{Cst}. To
study the Neimark-Sacker bifurcation, we let the unique interior
steady state be denoted by $E^*=(x^*, y^*)$.

\medskip

\noindent{\bf Theorem 4.3} {\em Let $\lambda>1$, $2\alpha\leq 1$
and $\beta \bar x>1$. Then the unique interior steady state
$E^*=(x^*, y^*)$ is asymptotically stable if $\beta<\beta_d$ and a
repeller if $\beta>\beta_d$. Moreover, $E^*$ undergoes a
Neimark-Sacker bifurcation at $\beta=\beta_d$. The bifurcation is
supercritical if $C^*>0$ and the bifurcation is subcritical if
$C^*<0$.}
\begin{proof} It remains to prove that a Neimark-Sacker
bifurcation  \cite{hale} occurs at $\beta=\beta_d$. Let
$\lambda_{\pm}$ denote the eigenvalues of $J(E^*)$ and $G$ be the
map induced by \eqref{rmodel}. We need to verify (a) $
 G(\beta, E^*)=E^*$ for $\beta$ near $\beta_d$, (b) $J(E^*)$ has two non-real eigenvalues for $\beta$ near $\beta_d$ with modulus $1$
 at $\beta=\beta_d$, (c) $\cfrac{d|\lambda_\pm|}{d\beta}>0$ at $\beta=\beta_d$, and (d) $\lambda_\pm^n\neq 1$ at $\beta=\beta_d$
for $n=1,2,3,4$.

It is clear that condition (a) holds. At $\beta=\beta_d$,
$det(J)=1$ and $|tr(J)|<1+det(J)=2$ imply
$\lambda_{\pm}=\cfrac{tr(J)\pm \sqrt{4-(tr(J))^2}i}{2}$ and
$|\lambda_{\pm}|=1$.
  As $\beta$ is varied around $\beta_d$, $J( E^*)$ varies continuously with respect to $\beta$. Hence eigenvalues of $J(E^*)$ are complex if
  $\beta$ is close to $\beta_d$ and condition (b) is true. To verify (c), notice
$\cfrac{d(det(J))}{d\beta}=\cfrac{d(det(J))}{dy^*}\times\cfrac{dy^*}{d\beta}>0$,
  and   hence  $\cfrac{d|\lambda_\pm|}{d\beta}\bigg|_{\beta=\beta_d}=\cfrac{1}{2}\cfrac{d(det(J))}{d\beta}\bigg|_{\beta=\beta_d}>0$.
 Therefore, eigenvalues of $J(E^*)$  cross the unit circle transversally. It
remains to verify $\lambda_{\pm}^n \neq 1$ for $n=1,2,3,4$.
Clearly, $\lambda_{\pm} \neq \pm 1$, and $\lambda_{\pm} = \pm i$
if and only if $tr(J)=0$,  which is impossible. Thus
$\lambda_{\pm}^n \neq 1$ for $n=1,2,4$. Also $\lambda_{\pm}^3=1$
if and only if  $tr(J)=-1$ at $\beta=\beta_d$ and we obtain a
contradiction.  Therefore, $\lambda_{\pm}^n \neq 1$ for $n=1,2,3,
4$ and a Neimark-Sacker bifurcation occurs at $\beta=\beta_d$ by
\cite{hale}.

To determine whether the Neimark-Sacker bifurcation is
supercritical or subcritical, we  perform a standard analysis as
we do in \cite{chow2016}. We first move the unique interior steady
state $(x^*, y^*)$ to the origin by letting $X=x-x^*$ and
$Y=y-y^*$, i.e.,
\begin{equation}\label{genori}\left\{\begin{array}{ll}
X_{t+1}= \cfrac{\lambda(X_t+x^*)}{1+X_t+x^*}\displaystyle e^{-(Y_t+y^*)(1+\alpha Y_t+\alpha y^*)}-x^*\\[2.5ex]
Y_{t+1}=\beta(X_t+x^*)\left(1-e^{-(Y_t+y^*)(1+\alpha Y_t+\alpha
y^*)}\right)-y^*.
\end{array}\right.
\end{equation}
Using the Taylor series expansion,  \eqref{genori} can be put
 into the following form
\begin{equation}\label{newNS}\left(\begin{array}{l}
 X_{t+1} \\
 \\
 Y_{t+1}
\end{array}\right)=J(E^*)
\left(\begin{array}{l}
 X_{t} \\
 \\
 Y_{t}
 \end{array}\right)+\left(\begin{array}{l}
 \hat{f}(X_t,Y_t) \\
 \\
 \hat{g}(X_t,Y_t)
 \end{array}\right),\end{equation}
where
\begin{equation}\begin{array}{l}
\hat{f}(X,Y)=b_1X^2+b_2XY+b_3Y^2+b_4X^3+b_5X^2Y+b_6XY^2+b_7Y^3+O(4)\\[1ex]
\hat{g}(X,Y)=c_1XY+c_2Y^2+c_3Y^3+c_4XY^2+O(4)
 \end{array}\end{equation}
with $$b_1=\cfrac{-\lambda e^{-\diamondsuit^*}}{(1+x^*)^3}, \
b_2=\cfrac{-\lambda (1+2\alpha
y^*)e^{-\diamondsuit^*}}{(1+x^*)^2}, \ b_3=\cfrac{\lambda
x^*e^{-\diamondsuit^*}[(1+2\alpha y^*)^2-2\alpha]}{2(1+x^*)},$$
$$ b_4=\cfrac{\lambda ^2e^{-\diamondsuit^*}}{(1+x^*)^4}, \
b_5=\cfrac{2\lambda (1+2\alpha
y^*)e^{-\diamondsuit^*}}{(1+x^*)^3}, \ b_6=\cfrac{\lambda
e^{-\diamondsuit^*}[(1+2\alpha y^*)^2-2\alpha]}{2(1+x^*)^2},$$
$$b_7=\cfrac{\lambda
x^*e^{-\diamondsuit^*}(1+2\alpha y^*)[-(1+2\alpha
y^*)^2+6\alpha]}{6(1+x^*)},$$ and
$$c_1=\beta e^{-\diamondsuit^*}(1+2\alpha y^*), \ c_2=\cfrac{-\beta
x^*e^{-\diamondsuit}[(1+2\alpha y^*)^2-2\alpha]}{2}$$
$$c_3=\cfrac{-\beta x^*e^{-\diamondsuit^*}(1+2\alpha
y^*)[-(1+2\alpha y^*)^2+6\alpha]}{6}, \ c_4=\cfrac{-\beta
e^{-\diamondsuit^*}[(1+2\alpha y^*)^2-2\alpha]}{2}.$$

Let $J(E^*)|_{\beta=\beta_d}=\left(\begin{array}{ll} a_{11} &
a_{12}\\
a_{21} & a_{22}\end{array}\right).$ Then $J(E^*)|_{\beta=\beta_d}$
has eigenvalues $\mu\pm i\omega$, where  $\omega>0$ and
$\mu^2+\omega^2=1$. Let $L=\left(\begin{array}{ll} a_{12} &
0\\
\mu-a_{11} & -\omega\end{array}\right)$ and define the new
variables $u$ and $v$ via
\begin{equation}\label{lin}\left(\begin{array}{l}
X_t\\Y_t\end{array}\right)=L\left(\begin{array}{l}u_t\\v_t\end{array}\right).\end{equation}
System \eqref{newNS} in terms of $u$ and $v$ becomes
\begin{equation}\left(\begin{array}{l}u_{t+1}\\v_{t+1}\end{array}\right)=\left(\begin{array}{lc}\mu
& -\omega\\\omega &
\mu\end{array}\right)\left(\begin{array}{l}u_t\\v_t\end{array}\right)+\left(\begin{array}{l}{\tilde f}(u_t,v_t)\\
{\tilde g}(u_t,v_t)\end{array}\right),
\end{equation}
where $${\tilde
f}=\cfrac{1}{a_{12}}(k_1u^2+k_2v^2+k_3uv+k_4u^3+k_5v^3+k_6u^2v+k_7uv^2)+O(4)$$
 $${\tilde
g}=l_1u^2+l_2v^2+l_3uv+l_4u^3+l_5v^3+l_6u^2v+l_7uv^2+O(4),$$ and
$$k_1=b_1a_{12}^2+b_2a_{12}(\mu-a_{11})+b_3(\mu-a_{11})^2, \
k_2=b_3\omega^2, k_3=-b_2a_{12}\omega+2b_3\omega(a_{11}-\mu),$$
$$k_4=b_4a_{12}^3+b_5a_{12}^2(\mu-a_{11})+b_6a_{12}(\mu-a_{11})^2+b_7(\mu-a_{11})^3,
k_5=-b_7\omega^3,$$
$$k_6=-b_5a_{12}^2\omega-2b_6a_{12}\omega(\mu-a_{11})-3b_7\omega(\mu-a_{11})^2, \ k_7=b_6a_{12}\omega^2+3b_7(\mu-a_{11})\omega^2,$$
$$l_1=\cfrac{\mu-a_{11}}{a_{12}\omega}k_1-\cfrac{1}{\omega}[c_1a_{12}(\mu-a_{11})+c_2(\mu-a_{11})^2], \ l_2=\cfrac{b_3(\mu-a_{11})}{a_{12}}\omega-c_2\omega,$$
$$l_3=\cfrac{\mu-a_{11}}{a_{12}\omega}k_3+c_1a_{12}+2c_2(\mu-a_{11}),
\
l_4=\cfrac{\mu-a_{11}}{a_{12}\omega}k_4-\cfrac{(\mu-a_{11})^2}{\omega}[c_3(\mu-a_{11})+c_4a_{12}],$$
$$l_5=-\cfrac{\mu-a_{11}}{a_{12}}b_7\omega^2+c_3\omega^2, \
l_6=\cfrac{\mu-a_{11}}{a_{12}\omega}k_6+(\mu-a_{11})[3(\mu-a_{11})+2c_4a_{12}],$$
$$l_7=\cfrac{\mu-a_{11}}{a_{12}\omega}k_7-[3c_3(\mu-a_{11})+c_4a_{12}]\omega$$

Applying Theorem 15.31 of \cite{hale}, the direction of a
Neimark-Sacker bifurcation is determined by
 \begin{equation}\label{Cst}
 C^*=Re\big(\cfrac{(1-2\lambda_+)\lambda_-^2}{1-\lambda_+}\xi_{20}\xi_{11}\big)+\cfrac{1}{2}|\xi_{11}|^2+|\xi_{02}|^2-Re(\lambda_-\xi_{21}),
 \end{equation}
 where $Re$ denotes the real part of a complex number and
\begin{eqnarray}\xi_{20}&=&\cfrac{1}{8}({\tilde f}_{uu}-{\tilde
f}_{vv}+2{\tilde
g}_{uv}+i({\tilde g}_{uu}-{\tilde g}_{vv}-2{\tilde f}_{uv}))\bigg|_{(0,0)}\nonumber\\
&=& \cfrac{1}{8}\big[(2k_1/a_{12}-2k_2/a_{12}+2l_3)+i(2l_1-2l_2-2k_3/a_{12})\big]\nonumber\\
 \xi_{11}&=&\cfrac{1}{4}({\tilde f}_{uu}+{\tilde f}_{vv}+i({\tilde
 g}_{uu}+{\tilde
 g}_{vv}))\bigg|_{(0,0)}\nonumber\\
 &=&
 \cfrac{1}{4}\big[(2k_1/a_{12}+2k_2/a_{12})+i(2l_1+2l_2)\big]\\
 \xi_{02}&=&\cfrac{1}{8}({\tilde f}_{uu}-{\tilde f}_{vv}-2{\tilde g}_{uv}+i({\tilde g}_{uu}-{\tilde
 g}_{vv}+2{\tilde
 f}_{uv}))\bigg|_{(0,0)}\nonumber\\
 &=&
 \cfrac{1}{8}\big[(2k_1/a_{12}-2k_2/a_{12}-2l_3)+i(2l_1-2l_2+2k_3/a_{12})\big]\nonumber\\
 \xi_{21}&=&\cfrac{1}{16}({\tilde f}_{uuu}+{\tilde f}_{uvv}+{\tilde g}_{uuv}+{\tilde g}_{vvv}+i({\tilde g}_{uuu}+{\tilde g}_{vvv}-{\tilde
 f}_{uuv}-{\tilde
 f}_{vvv}))\bigg|_{(0,0)}\nonumber\\
 &=&
 \cfrac{1}{16}\big[(6k_4/a_{12}+2k_7/a_{12}+2l_6+6l_5)+i(6l_4+6l_5-2k_6/a_{12}-6k_5/a_{12})\big]\nonumber.
 \end{eqnarray}
 If $C^*>0$, then the system has an attracting closed invariant circle for $\beta>\beta_d$ and near $\beta_d$. If $C^*<0$,
 then the system has an unstable closed invariant circle for $
 \beta<\beta_d$ and near $\beta_d$. \end{proof}

\medskip
It is not easy to determine analytically whether the bifurcation
is supercritical or subcritical since $C^*$ cannot be computed
analytically. Numerical investigation does indicate that the
bifurcation is supercritical so that the model has an attracting
invariant closed curve when $\beta>\beta_d$ and near $\beta_d$.

We conclude from Theorems 4.2, 4.3 and Proposition 2.3 that
cooperative hunting does not affect dynamical interactions of the
prey and predator if the degree of cooperative hunting is small,
$2\alpha\leq 1$.

\subsection{Dynamics of the model when $2\alpha>1$}
Recall Theorem 3.1 indicates that \eqref{rmodel} has a unique
interior steady state if $\beta \bar x>1$, \eqref{rmodel} has no
interior steady state if $\beta \bar x<1$ and $1<2\alpha\leq
\cfrac{3\lambda-1}{\lambda-1}$, and there are either zero, one or
two interior steady states if $\beta \bar x<1$ and
$2\alpha>\cfrac{3\lambda-1}{\lambda-1}$. We shall determine
stability of an interior steady state when it exists. In the case
when \eqref{rmodel} has no interior steady state, we suspect that
$E_1$ is globally asymptotically stable. The following result
provides a restriction on the parameter $\alpha$ for which $E_1$
is  a global attractor.

\medskip

\noindent{\bf Theorem 4.4} {\em Let $\lambda>1$, $2\alpha>1$ and
$\beta \bar x \sqrt{2\alpha}e^{\cfrac{1-2\alpha}{4\alpha}}<1$.
Then $E_1=(\bar x, 0)$ is globally asymptotically stable in
$\{(x,y)\in\mathbb{R}_+^2: x>0\}$.}

\begin{proof} We first observe that
$2\alpha e^{\cfrac{1-2\alpha}{2\alpha}}>1$. Indeed, since
$e^{1-1/x}<x$ for $x>1$, letting $x=2\alpha$, we obtain the
inequality. Using the assumption,  $\beta \bar x<\beta \bar x
\sqrt{2\alpha}e^{\cfrac{1-2\alpha}{4\alpha}}<1$  and $E_1$ is
locally asymptotically stable by Proposition 2.3.

Let $(x(0),y(0))$ be given with $x(0)>0$. We may assume $y(0)>0$
and thus $y(n)>0$ for $n\geq 0$. Since
$\limsup_{n\rightarrow\infty}x(n)\leq \bar x$, for any
$\epsilon>0$ there exists $n_0>0$ such that $x(n)<\bar x+\epsilon$
for $n\geq n_0$. We choose $\epsilon>0$ so that $\beta (\bar
x+\epsilon) \sqrt{2\alpha}e^{\cfrac{1-2\alpha}{4\alpha}}<1$. Then
from the second equation of \eqref{rmodel}, we have
$y(n+1)<\beta(\bar x+\epsilon)\bigg(1-e^{-y(n)(1+\alpha
y(n))}\bigg)$ for $n\geq n_0$. Consider the scalar equation
\begin{equation}
z(n+1)=\beta(\bar x+\epsilon)\bigg(1-e^{-z(n)(1+\alpha
z(n))}\bigg), \ z(0)=y(n_0)
\end{equation}
and letting $l(z)=\beta(\bar x+\epsilon)(1-e^{-z(1+\alpha z)})$.
It can be easily shown that $l(z)<z$ for $z>0$. Indeed, letting
$R(z)=l(z)-z$. Then $R(0)=0$, $R^\prime(z)=\beta (\bar
x+\epsilon)e^{-z(1+\alpha z)}(1+2\alpha z)-1$ with
$R^\prime(0)=\beta (\bar x+\epsilon)-1<0$ and
$R^{\prime\prime}(z)=\beta (\bar x+\epsilon)e^{-z(1+\alpha
z)}\bigg(2\alpha-(1+2\alpha z)^2\bigg)$. Now
$R^{\prime\prime}(z)=0$ has a unique positive solution
$s^*=\cfrac{1}{2\alpha}(\sqrt{2\alpha}-1)$ such that
$R^{\prime\prime}(z)>0$ on $(0, s^*)$ and $R^{\prime\prime}(z)<0$
 on $(s^*, \infty)$. Therefore, the maximum value
of $R^\prime(z)$ attains at $z=s^*$.  A simple calculation yields
$R^\prime(s^*)=\beta  (\bar
x+\epsilon)\sqrt{2\alpha}e^{\cfrac{1-2\alpha}{4\alpha}}-1<0$ by
the assumption. Therefore, $R^\prime(z)<0$ for $z>0$ and $l(z)<z$
is shown. It follows that $\lim_{n\rightarrow\infty}z(n)=0$ if
$z(0)\geq 0$ and hence $\lim_{n\rightarrow\infty}y(n)=0$ for
$y(0)\geq 0$. Consequently, $\lim_{n\rightarrow\infty}x(n)=\bar x$
if $x(0)>0$ and the proof is complete. \end{proof}

Theorem 4.4 provides a sufficient condition for which the
predators go extinct when the degree of cooperation is large. That
is, large cooperation among predators drive the predators to
extinction under the condition given by the theorem.  We now study
local stability of an interior steady state $E=(x,y)$ when
$1<2\alpha$. Recall from Lemma 4.1 and \eqref{Jury} that in order
to determine local stability at $E$, we have to study $V(y)$
defined in \eqref{V}. Note that $y_c$ is defined in \eqref{y_c}.
The result is  as follows.

\medskip

\noindent{\bf Lemma 4.5} {\em Let $\lambda>1$ and $2\alpha>1$.
Then $V(0)=0$. Moreover, $V(y)>0$ for $y\in(0, y_c)$ if
$2\alpha\leq\cfrac{3\lambda-1}{\lambda-1}$ and there exists a
unique $y_t\in(0, y_c)$ such that $V(y)<0$ on $(0, y_t)$ and
$V(y)>0$ on $(y_t, y_c)$ if
$2\alpha>\cfrac{3\lambda-1}{\lambda-1}$.}
\begin{proof}
 Using \eqref{J1st}, a direct computation yields
\begin{equation}\label{Vex}
V(y)= 1-\cfrac{e^\diamondsuit}{\lambda} +\bigg \{
2-(1-\cfrac{1}{\lambda})\cfrac{1}{1-e^{-\diamondsuit}}\bigg
\}y(1+2\alpha y).
\end{equation}
It is easy to check that $V(0)=0$ and
$V^\prime(0)=\cfrac{(3\lambda-1)-2\alpha(\lambda-1)}{2\lambda}$.
We need to solve $V(y)=0$ on $y \in (0, y_c)$.  Notice
$1-\cfrac{e^\diamondsuit}{\lambda} >0$ on  $(0, y_c)$. Let $\tilde
y_1 \in (0, y_c) $ satisfying
$e^\diamondsuit=\cfrac{2\lambda}{\lambda+1}$. Thus $e^\diamondsuit
\geq \cfrac{2\lambda}{\lambda+1}$ for $ y \geq \tilde y_1$. We
first show that
\begin{equation}\label{Vpos}
V(y)>0 \mbox{ on } y \in [\tilde y_1, y_c]
\end{equation}
as  the term inside the bracket of \eqref{Vex} equals $   \cfrac{(\lambda+1)e^{\diamondsuit}-2 \lambda }{\lambda ( e^{\diamondsuit}-1)}  $
which is positive on $(\tilde y_1, y_c]$ and negative on $(0, \tilde y_1)$.

By  \eqref{Vpos} we need only to consider  $y \in (0, \tilde y_1)$
as follows. Define $F(y)  =\cfrac{(
e^{\diamondsuit}-1)(\lambda-e^{\diamondsuit}) }{2 \lambda
-(\lambda+1)e^{\diamondsuit}}-y(1+2\alpha y).$ Then
\begin{equation}\label{F}
F(y)  = \cfrac{\lambda ( e^{\diamondsuit}-1)}{2 \lambda
-(\lambda+1)e^{\diamondsuit} } V(y),  \mbox{  and thus sign }V
=\mbox{ sign }F \mbox{ for } y \in (0, \tilde y_1).
\end{equation}
Moreover,
\begin{equation}\label{F0}
F(0)=0 \mbox{ and } F(\tilde y_1^{-}) = \infty.
\end{equation}
Differentiating $F(y)$,
$$
F^\prime(y)=\cfrac{(\lambda+1)e^{3\diamondsuit}-4\lambda
e^{2\diamondsuit}+\lambda(\lambda+1)e^{\diamondsuit}}{[2\lambda-(\lambda+1)e^{\diamondsuit}]^2}(1+2\alpha
y)-(1+4\alpha y)
$$
with $F^{\prime}(0)=0$ and $F^{\prime}(\tilde y_1^{-}) = \infty$.
Define $G(y)= \cfrac{F^{\prime}(y)}{a(1+2\alpha y)}$, i.e.,
\begin{equation}\label{G}
G(y)  =\cfrac{(\lambda+1)e^{3\diamondsuit}-4\lambda
e^{2\diamondsuit}+\lambda(\lambda+1)e^{\diamondsuit}}{[2\lambda-(\lambda+1)e^{\diamondsuit}]^2}-\cfrac{1+4\alpha
y}{1+2\alpha y}.
\end{equation}
Similar to  \eqref{F},
\begin{equation}\label{signG}
\mbox{sign }F^{\prime} =\mbox{ sign }G   \mbox{ for } y \in (0, \tilde y_1)
\end{equation}
with $G(0)=0$ and $G(\tilde y_1^{-}) = \infty$. A simple calculation shows
\begin{equation}\label{L}
\cfrac{G^\prime(y)}{(1+2\alpha
y)}=\cfrac{e^{\diamondsuit}}{\lambda+1}\bigg(1+\lambda(\lambda-1)^2\cfrac{(\lambda+1)e^{\diamondsuit}+2\lambda}
{[2\lambda-(\lambda+1)e^{\diamondsuit}]^3}\bigg)-\cfrac{2\alpha}{(1+2\alpha
y)^3}.
\end{equation}
Denote by $L(y)$ the function on the right-hand side of \eqref{L}.
Then
\begin{equation}\label{L0}
L(0)=(\cfrac{3\lambda-1}{\lambda-1}-2\alpha) \geq 0 \mbox{ if and
only if }  2\alpha   \leq \cfrac{3\lambda-1}{\lambda-1}.
\end{equation}
We claim that
\begin{equation}\label{Lup}
L(\tilde y_1^{-}) = \infty \mbox{ and } L(y) \uparrow \mbox{ strictly on } (0, \tilde y_1).
\end{equation}
The first claim is easily seen from \eqref{L}. Notice that both
$\cfrac{e^{\diamondsuit}}{\lambda+1}$ and
$\cfrac{-2\alpha}{(1+2\alpha y)^3}$ are strictly increasing in
$y$. Let  $s=(\lambda+1)e^\diamondsuit$, which is strictly
increasing in $y$. It suffices to show that
$$
\cfrac{(\lambda+1)e^{\diamondsuit}+2\lambda}{[2\lambda-(\lambda+1)e^{\diamondsuit}]^3}=\cfrac{s+2\lambda}{(2\lambda-s)^3}  \mbox{ is increasing for } 0<s< 2 \lambda.
$$
Note that $0<s< 2 \lambda$ is equivalent to $y<\tilde y_1$ as $\tilde y_1$ is defined by
$e^\diamondsuit=\cfrac{2\lambda}{\lambda+1}$.
Since   $\cfrac{d}{ds}   \cfrac{s+2\lambda}{(2\lambda-s)^3}  =\cfrac{2s+8\lambda}{(2\lambda-s)^4}>0$,  \eqref{Lup} is verified.

If  $2\alpha   \leq \cfrac{3\lambda-1}{\lambda-1} $,  \eqref{L0}
and \eqref{Lup} imply $L(y) >0$ on $(0, \tilde y_1)$. Then using
\eqref{L} backward one by one up to \eqref{F}, we obtain
$G^{\prime}>0, G>0,F^{\prime}>0, F>0$ and at last $V>0$ on  $(0,
\tilde y_1)$. Together with \eqref{Vpos},  that $V(y)>0$ for
$y\in(0, y_c)$ is verified if
$2\alpha\leq\cfrac{3\lambda-1}{\lambda-1}$.

It remains to consider the case  $2\alpha  >
\cfrac{3\lambda-1}{\lambda-1}$. By \eqref{L0}  and \eqref{Lup},
$L(y) <0$ on $(0, \delta)$ and
 $L(y) >0$ on $(\delta, \tilde y_1)$ for some $\delta \in  (0, \tilde y_1)$. The same holds for $G^{\prime}$ by \eqref{L}.
 Since $G(0)=F^{\prime}(0)=F(0) =0$  and $G(\tilde y_1^{-}) = F^{\prime}(\tilde y_1^{-}) =  F(\tilde y_1^{-}) = \infty$,
 the same result on $(0, \tilde y_1)$ holds for $ G, F^{\prime} $ and then $F $ except probably each with a different constant $\delta$.
The conclusion for the present case follows from \eqref{F} and
\eqref{Vpos}.
 \end{proof}

\medskip
Figure 2(b)  plots $V(y)$ using $\lambda=10$, $\beta=6.3/20$, and
$\alpha=1.5/2.1$. Then
$2\alpha-\cfrac{3\lambda-1}{\lambda-1}=-1.7937<0$ and hence
$V(y)>0$ on $(0, y_c)$. In Figure 2(c), $\alpha$ is increased to
$15/2.1$ so that
$2\alpha-\cfrac{3\lambda-1}{\lambda-1}=11.0635>0$, and thus
$V(y)=0$ has a solution $y_t$ in $(0, y_c)$.

  If
$1<2\alpha\leq\cfrac{3\lambda-1}{\lambda-1}$ and $\beta  \bar
x>1$, then \eqref{rmodel} has a unique interior steady state
$E^*=(x^*, y^*)$ by Theorem 3.1. It follows from Lemma 4.5 that
$tr(J)<1+det(J)$ and there exists a unique $y_d>0$, $y_d<y_c$,
such that $det(J)|_{y=y_d}=1$ by Lemma 4.1. Since the $y$
component of any interior steady state is a strictly increasing
function of $\beta$, there exists a unique $\beta_d>0$ such that
$E^*$ is asymptotically stable if $\beta\in(0, \beta_d)$ and a
repeller if $\beta>\beta_d$. Similar to the proof of Theorem 4.2,
it can be easily verified that $E^*$ undergoes a Neimark-Sacker
bifurcation at $\beta=\beta_d$. We summarize the discussion as
follows.

\medskip

\noindent{\bf Theorem 4.6} {\em Let $\lambda>1$, $\beta  \bar x>1$
and  $1<2\alpha\leq\cfrac{3\lambda-1}{\lambda-1}$. Then
\eqref{rmodel} has a unique interior steady state $E^*=(x^*,
y^*)$, where $E^*$ is asymptotically stable if $\beta<\beta_d$ and
a repeller if $\beta>\beta_d$. Moreover, $E^*$ undergoes a
Neimark-Sacker bifurcation at $\beta=\beta_d$.}

\medskip

When the maximal reproductive number of the predators exceeds one
 and the degree of cooperation $\alpha$ is neither too
small nor too large, Theorem 4.6 implies that the predator-prey
interaction can support a unique coexisting steady state.
Consequently,  both populations can coexist indefinitely as a
steady state if the predator conversion $\beta$ is smaller than a
critical value $\beta_d$. Otherwise, coexistence  of both
populations may be more complicated if $\beta$ is larger than
$\beta_d$.

 Let $\lambda>1$ and
$2\alpha>\cfrac{3\lambda-1}{\lambda-1}$. Then \eqref{rmodel} has a
unique interior steady state if  $\beta  \bar x>1$ and the number
of interior steady states is either zero, one or two if $\beta
\bar x<1$. By Theorem 3.2, there exists a unique $\beta_*$ such
that there is no interior steady state if $\beta<\beta_*$. The
system has two interior steady states if
$\beta_*<\beta<\cfrac{1}{\bar x}$ and there is a unique interior
steady state if $\beta\geq\cfrac{1}{\bar x}$. See Fig 1(a).





We present some numerical investigations for the asymptotic
dynamics of the system. Using the parameter values $\lambda=5$ and
$\alpha=1/2.1$, then a unique interior steady state exists if
$\beta>0.25$ and a Neimark-Sacker bifurcation occurs when $\beta$
is close to $0.6$. Figure 3(a) presents an invariant closed curve
for $\beta=0.609$. We next increase $\alpha$ to $3/2.1$, Figure
3(b) provides an invariant closed curve for the case when $\beta
\bar x>1$ and $2\alpha<\cfrac{3\lambda-1}{\lambda-1}$. The initial
conditions are chosen near the  unstable unique interior steady
state for both plots.

We next study the scenario when
$2\alpha>\cfrac{3\lambda-1}{\lambda-1}$, $\beta \bar x<1$ and the
system has two interior steady states. The parameter values used
are $\lambda=5$, $\beta=4.2/20$, and $\alpha=20/2.1$. We choose an
initial condition $(2.3, 0.2)$ which is close to $E_2^*$. The
solution converges to a closed invariant circle. If initial
condition $(3.9, 0.1)$ is used, then the solution converges to the
boundary steady state $E_1=(\bar x, 0)=(4, 0)$ as shown in Fig
3(c). Notice that in this parameter regime, both $E_1^*$ and
$E_2^*$ are unstable. We wish to demonstrate the stability of
$E_2^*$ and thus we decrease $\beta$ to $3.76/20$ while keep all
other parameter values the same. Then the two isoclines have two
positive intersections which results in two interior steady
states. We use the same initial conditions as in Fig 3(c). In this
circumstance,   one solution converges to the stable interior
steady state while the other solution converges to the boundary
steady state $E_1$. See Fig 3(d). Therefore, bistability occurs
and the predator may survive depending on initial conditions while
the predator would go extinct if there is no cooperative hunting.

\medskip

We summarize conditions for the existence of interior steady
states in Table 1 and a list of notations is given in Table 2.

\bigskip

\noindent{\bf Table 1} Existence of interior steady states

\medskip

\begin{tabular}{ccc}\hline
 Parameter regime & parameter regime & number of interior steady states\\
 \hline
 $\beta \bar x>1$ &$\alpha\geq 0$ & 1\\[1ex]
 $\beta  \bar x<1$ & $2\alpha\leq 1$ & 0\\[1ex]
 $\beta  \bar x<1$ & $1<2\alpha\leq \cfrac{3\lambda-1}{\lambda-1}$ &
 0\\[1.5ex]
$\beta \bar x<1$ & $2\alpha>\cfrac{3\lambda-1}{\lambda-1}$ & 0, 1
or
 2\\
 \hline
 \end{tabular}


\bigskip

\noindent{\bf Table 2} List of notations

\medskip

\begin{tabular}{cc}
 Notation & definition\\
 \hline\
 $\bar x$ & $\lambda-1$ \\[1.5ex]
 $y_c$ & $e^{y_c(1+\alpha y_c)}=\lambda$\\[1.5ex]
 $y_d$ & $det(J)|_{y=y_d}=1$ \\[1.5ex]
 $y_t$ & $tr(J)|_{y=y_t}=1+det(J)|_{y=y_t}$\\[1.5ex]
 $\beta_d$ & $det(J)|_{\beta=\beta_d}=1$ \\[1.5ex]
$\diamondsuit$ & $y(1+\alpha y)$ \\[1.5ex]
$\tilde y_1$ & $e^{\diamondsuit}|_{y=\tilde y_1}=\cfrac{2\lambda}{\lambda+1}$\\
 \hline
 \end{tabular}

\bigskip

\noindent {\bf Remark.} In this investigation, we have not studied
stability of the interior steady states for
$2\alpha>\cfrac{3\lambda-1}{\lambda-1}$.
 The stability of such a steady state $E=(x,y)$ depends on the location of $y$ relative to $y_d$ and $y_t$. In addition, there are also $y_*$ and $y_e$ involved.
 See \eqref{order}, Lemmas 4.1 and 4.5.
 It is hard to compare the order of these quantities theoretically. We postpone our investigation to a future study.

\section{Summary and conclusions}
\setcounter{equation}{0} Mathematical models of predator-prey
interactions are interesting dynamical systems. There are many
populations in nature with non-overlapping generations.
Consequently, continuous-time models are not appropriate to
describe such population interactions and discrete-time systems
can be used to explore such populations.

Cooperation among individuals of the same predator species is
frequently observed in nature and it can change dynamical
interactions of biological systems \cite{scheel, uetz}. Motivated
by the recent research of Alves and Hilker \cite{alves} on
continuous-time models of predator-prey interactions with
cooperative hunting in  predators, we propose and investigate a
parallel discrete-time system. The model derivation is based on
the classical Nicholson-Bailey system but with density-dependent
growth rate in the prey population. Similar to \cite{alves},
cooperative hunting of the predator is modeled via the attack rate
of the predator. Due to this cooperation, the probability of an
individual prey escaped from being preyed upon is decreased. In
order to investigate the effects of cooperative hunting, the
dynamics of the system with no cooperation among predators are
summarized first.

Comparing the system of  cooperation with that of no cooperation,
several similar dynamical results are obtained. Indeed, both
populations go extinct if the intrinsic growth rate of prey is
smaller than one while both populations can coexist under the same
sufficient conditions, namely that the prey's intrinsic growth
rate and the predator's  maximal reproductive number are both
greater than. Further, it is proven that asymptotic dynamics of
the model are similar to the system with no cooperation if
$2\alpha\leq \cfrac{3\lambda-1}{\lambda-1}$, where $\alpha$ is the
degree of cooperation. Consequently,  if the degree of cooperation
$\alpha$ is small, then cooperative hunting does not change
dynamical interactions between the prey and predators. On the
other hand, if
 the degree of cooperation $\alpha$ is large, i.e., $2\alpha>\cfrac{3\lambda-1}{\lambda-1}$,
then cooperative hunting becomes critical for the survival of the
predator in the case that $\beta  \bar x<1$. The lumped parameter
$\beta  \bar x$ can be interpreted as the maximal reproductive
number of the predator. Without cooperation, the predator
population goes extinct if this reproductive number is smaller
than one. See Proposition 2.3. With cooperative hunting, the
predator-prey interactions may support two interior steady states
when this reproductive number is less than one as illustrated in
Theorem 3.1. As a result, the predator and prey may coexist even
if the maximal reproductive number of predator is smaller than
one. Therefore, cooperation between predators can promote survival
of the predator which would otherwise go  extinct in the absence
of this mechanism.

Comparing our results with those of the continuous-time model
studied by Alves and Hilker \cite{alves}, first notice that in the
absence of predator's cooperation the unique interior steady state
in \cite{alves} is globally asymptotically stable whenever it
exists. This is not true for our system since the unique interior
steady state  can undergo a Neimark-Sacker bifurcation  for system
\eqref{rmodel} when  predators do not engage in cooperation.
However, using the concept of uniform persistence we prove that
both populations can coexist indefinitely as long as the maximal
reproductive of predators exceed one independent of whether
predators cooperate or not. This coexistence is not proved in
\cite{alves} when predators engage in hunting cooperation. On the
other hand, the number of interior steady states for both of the
continuous and discrete-time models is the same. In particular,
both systems can have two coexisting steady state if the degree of
cooperation is large and  the predator's maximal reproductive
number of predator is less than one. In our study, however, we are
able to quantify this degree of cooperation explicitly in terms of
the prey's intrinsic growth rate. Furthermore, for small degree of
cooperation, the asymptotic dynamics of the continuous-time model
are the same as the model with no cooperation. The discrete-time
model proposed in this study also possess this property, namely
that the asymptotic dynamics of the system with small magnitude of
predator cooperation behave asymptotically the same as the model
of no cooperation.


\begin{figure}
\begin{center}
\includegraphics[height=2.0in, width=2.5in]{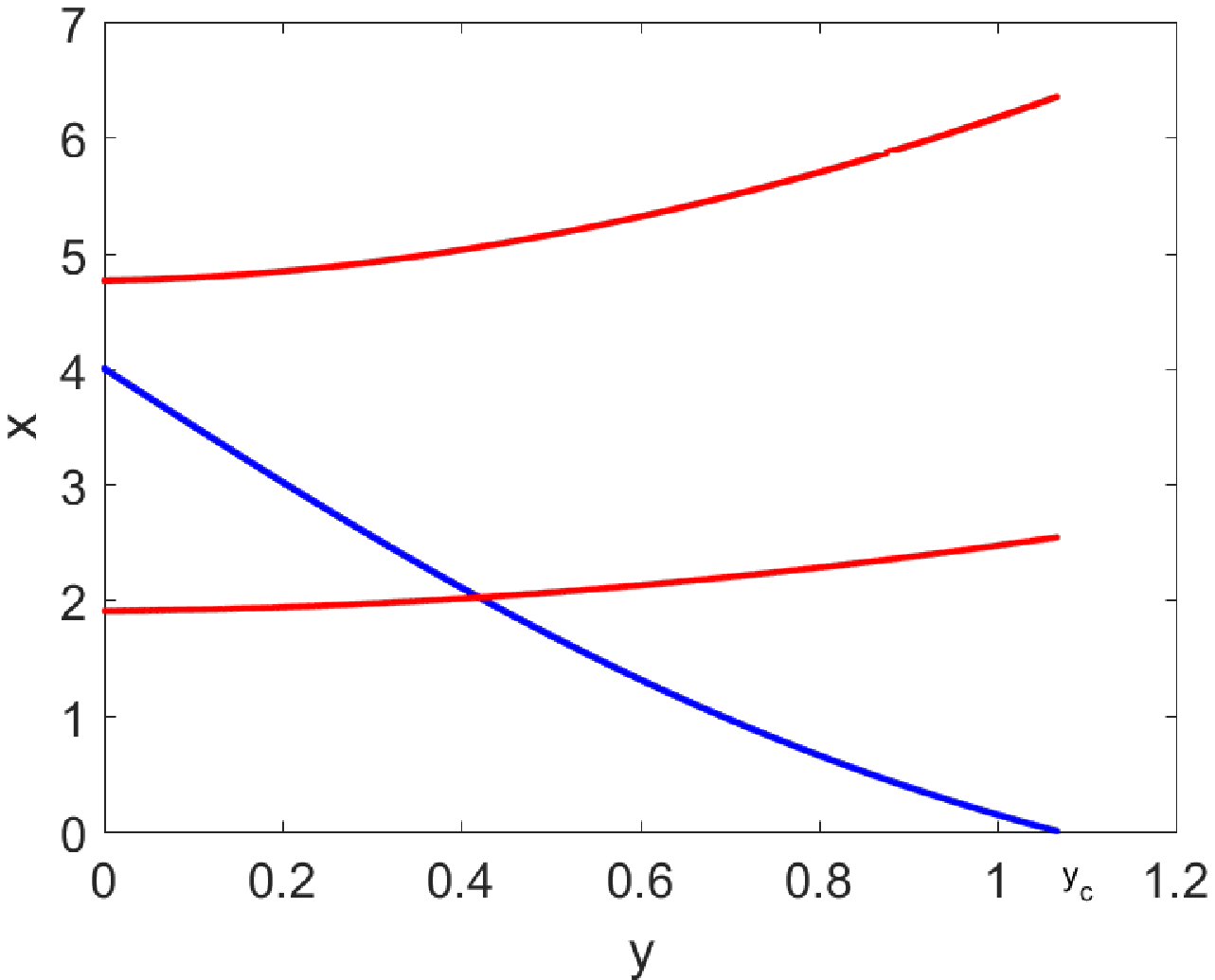} \\
(a)
\\
\begin{tabular}{cc}
\includegraphics[height=2.0in, width=2.5in]{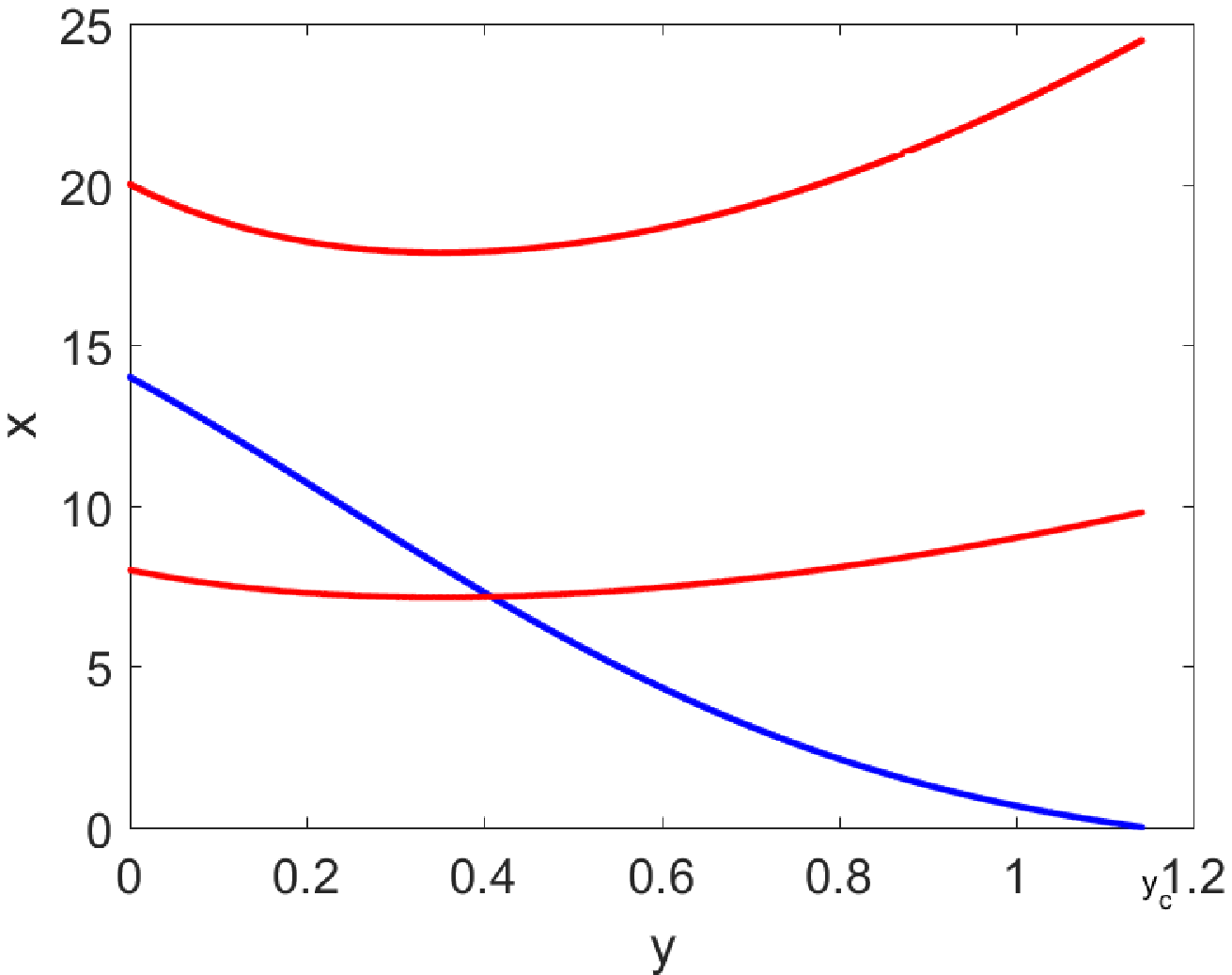}
&
\includegraphics[height=2.0in, width=2.5in]{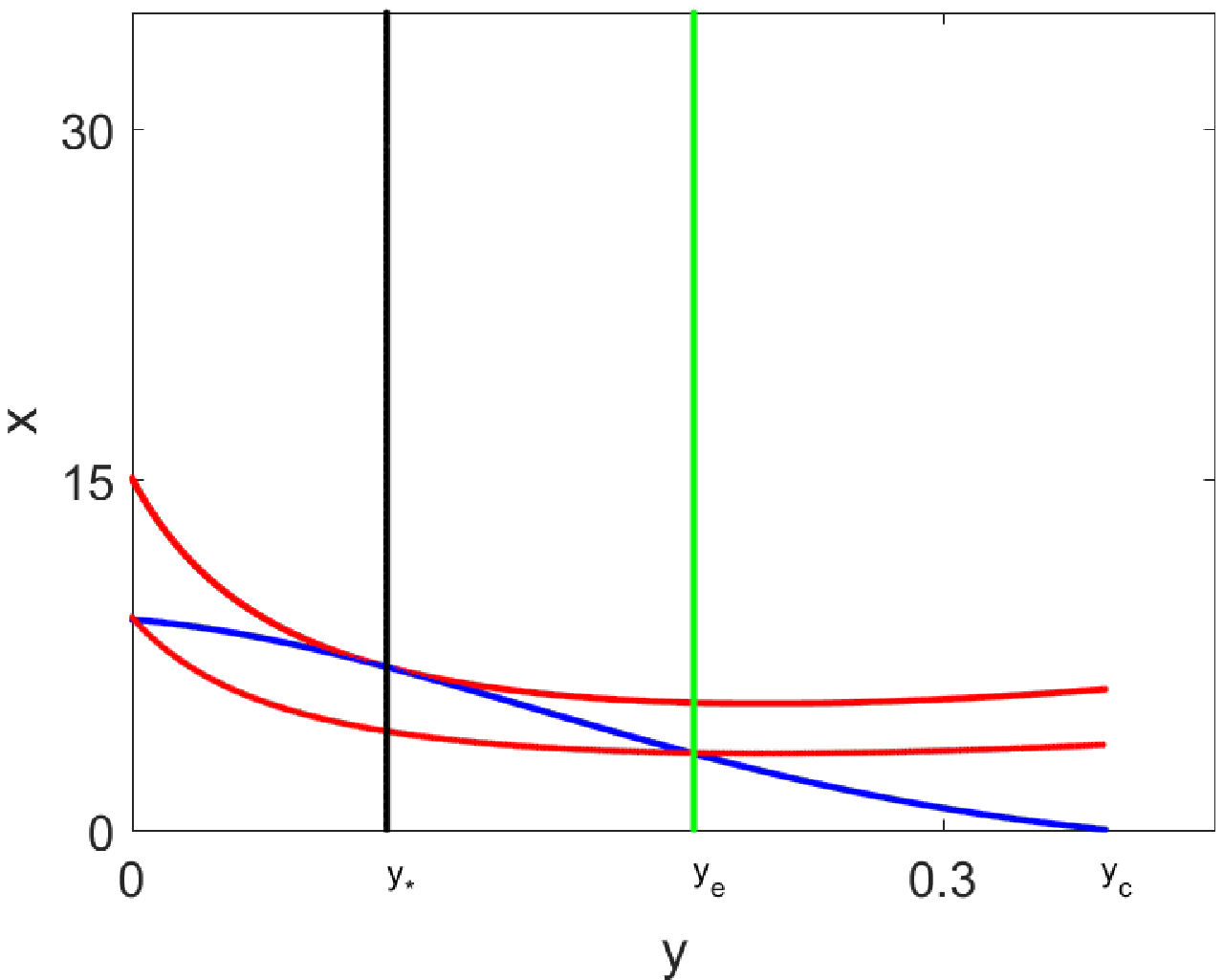}
\\
(b) & (c)
\end{tabular}
\end{center}

\caption{Isoclines are plotted for different $\beta$ values. In
(a) $2\alpha\leq 1$ while $2\alpha>1$ in (b) and (c). Further,
$2\alpha<\cfrac{3\lambda-1}{\lambda-1}$ in (b) and
$2\alpha>\cfrac{3\lambda-1}{\lambda-1}$ in (c).}
\end{figure}

\begin{figure}
\begin{center}
\includegraphics[height=2.0in, width=2.5in]{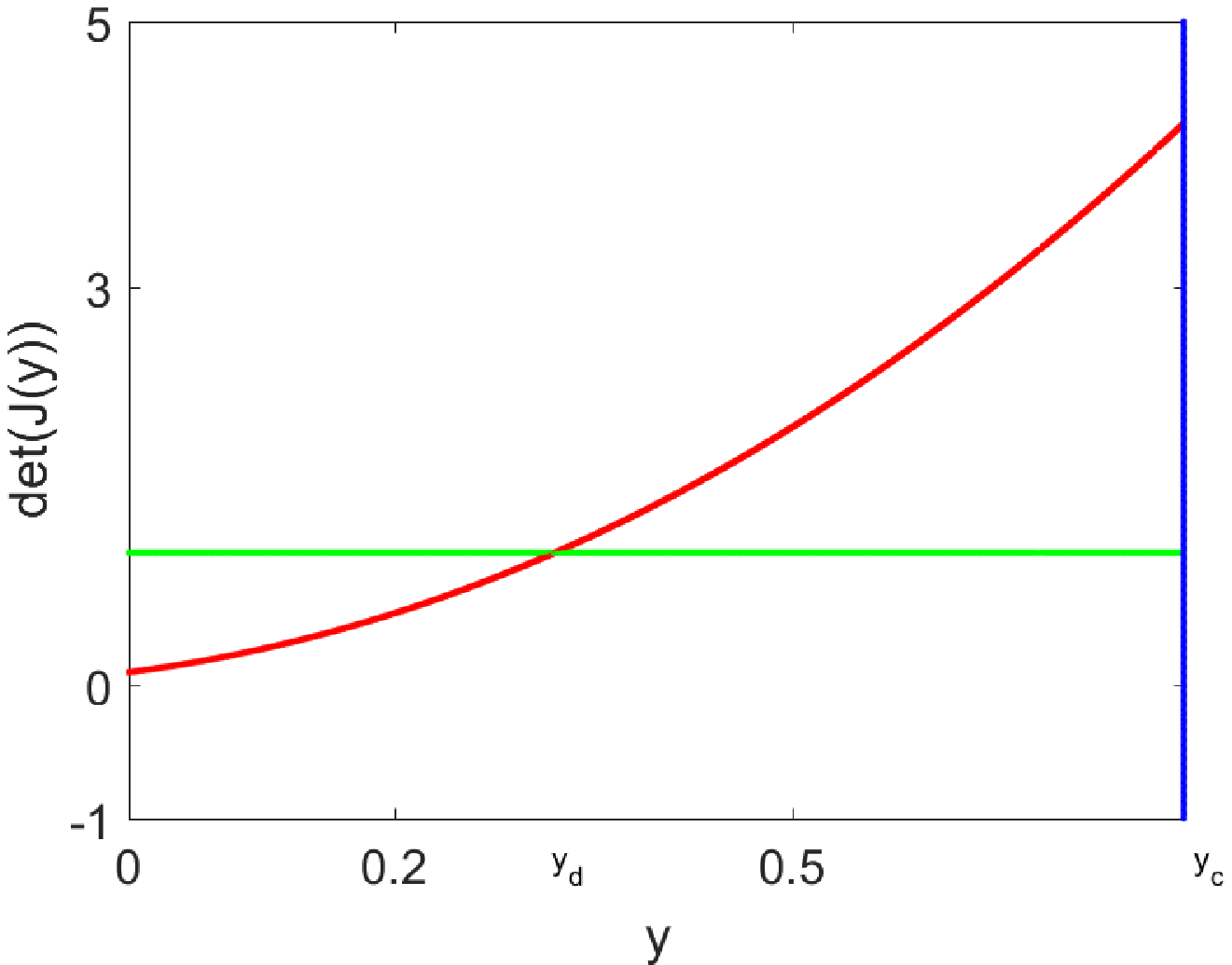}
\\
(a)
\\
\begin{tabular}{cc}
\includegraphics[height=2.0in, width=2.5in]{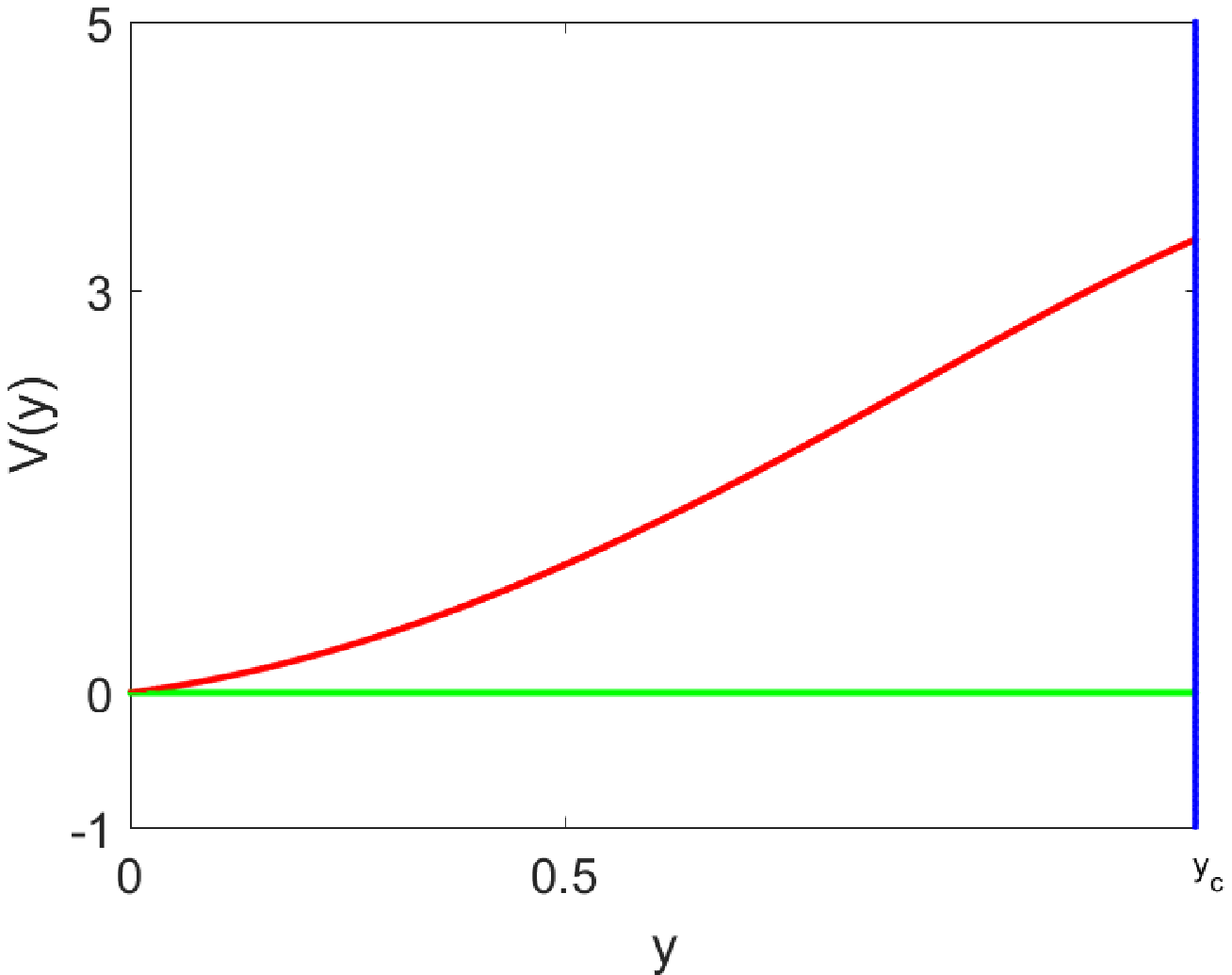} &
\includegraphics[height=2.0in, width=2.5in]{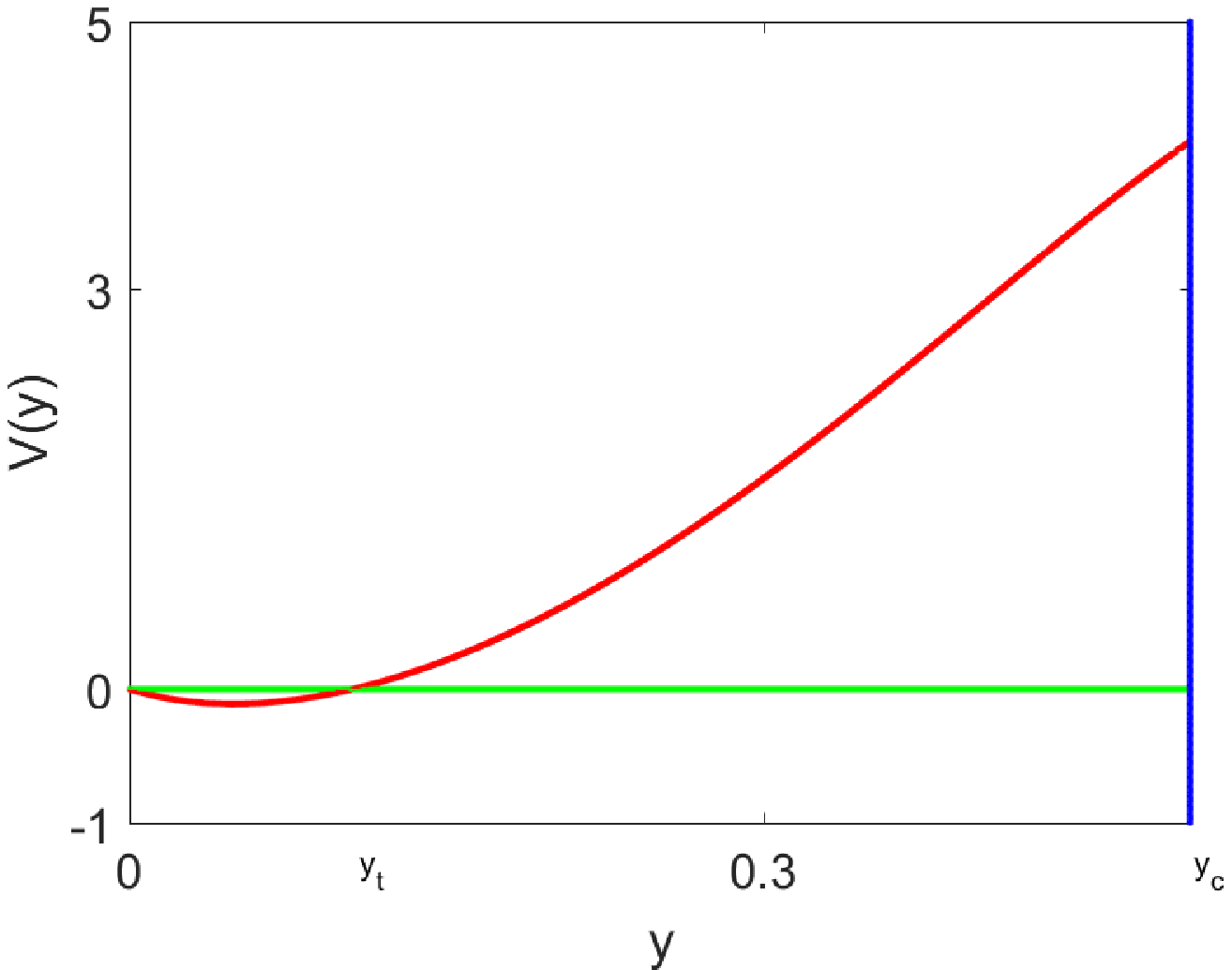}
\\
(b) & (c)
\end{tabular}
\end{center}
\caption{(a) plots $det(J)$ as a function of $y$  while (b) and
(c) provide the graphs of $V(y)$. In (b)
$2\alpha<\cfrac{3\lambda-1}{\lambda-1}$ and in (c)
$2\alpha>\cfrac{3\lambda-1}{\lambda-1}$.}
\end{figure}

\begin{figure}
\begin{center}
\begin{tabular}{cc}
\includegraphics[height=2.0in, width=2.5in]{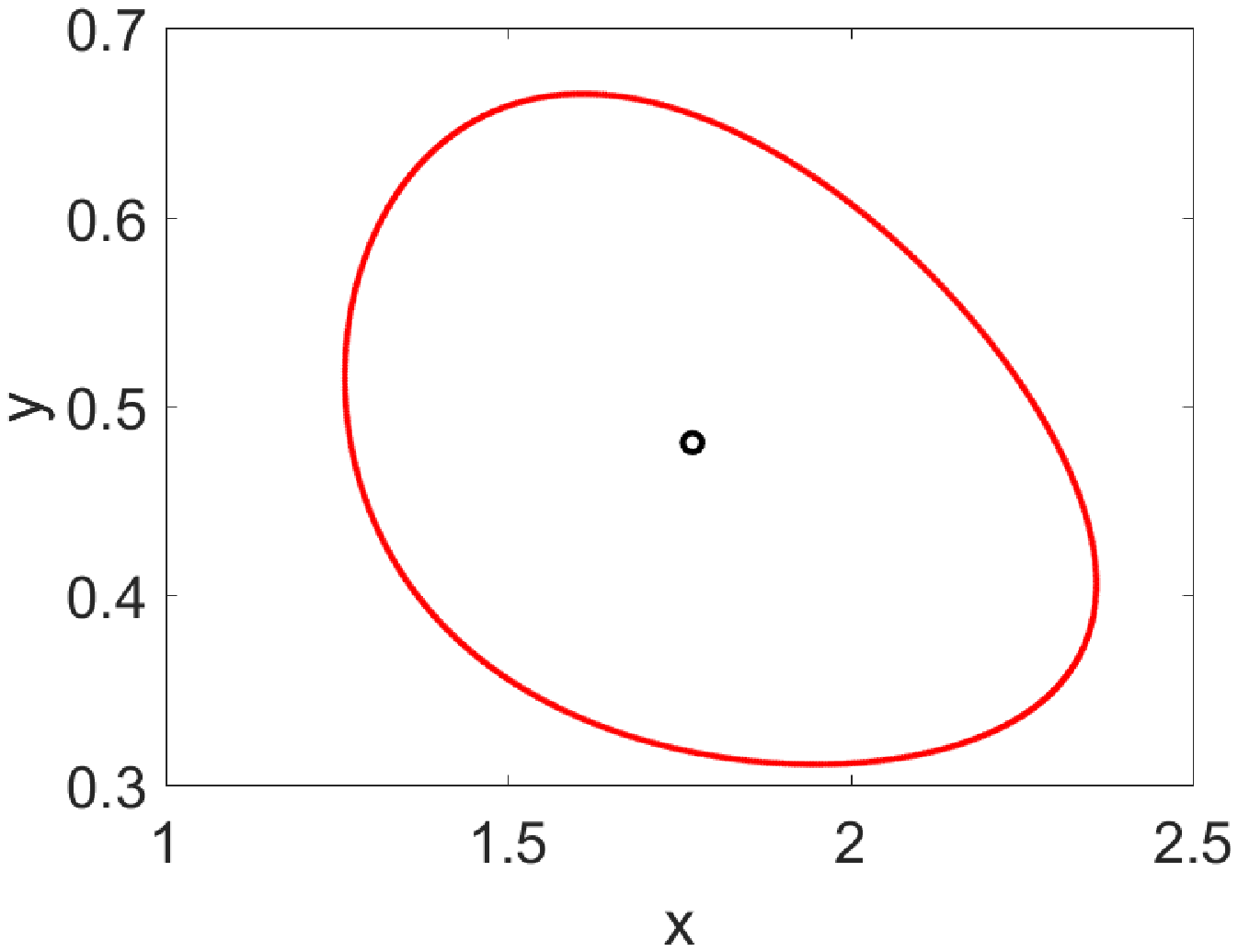} &
\includegraphics[height=2.0in, width=2.5in]{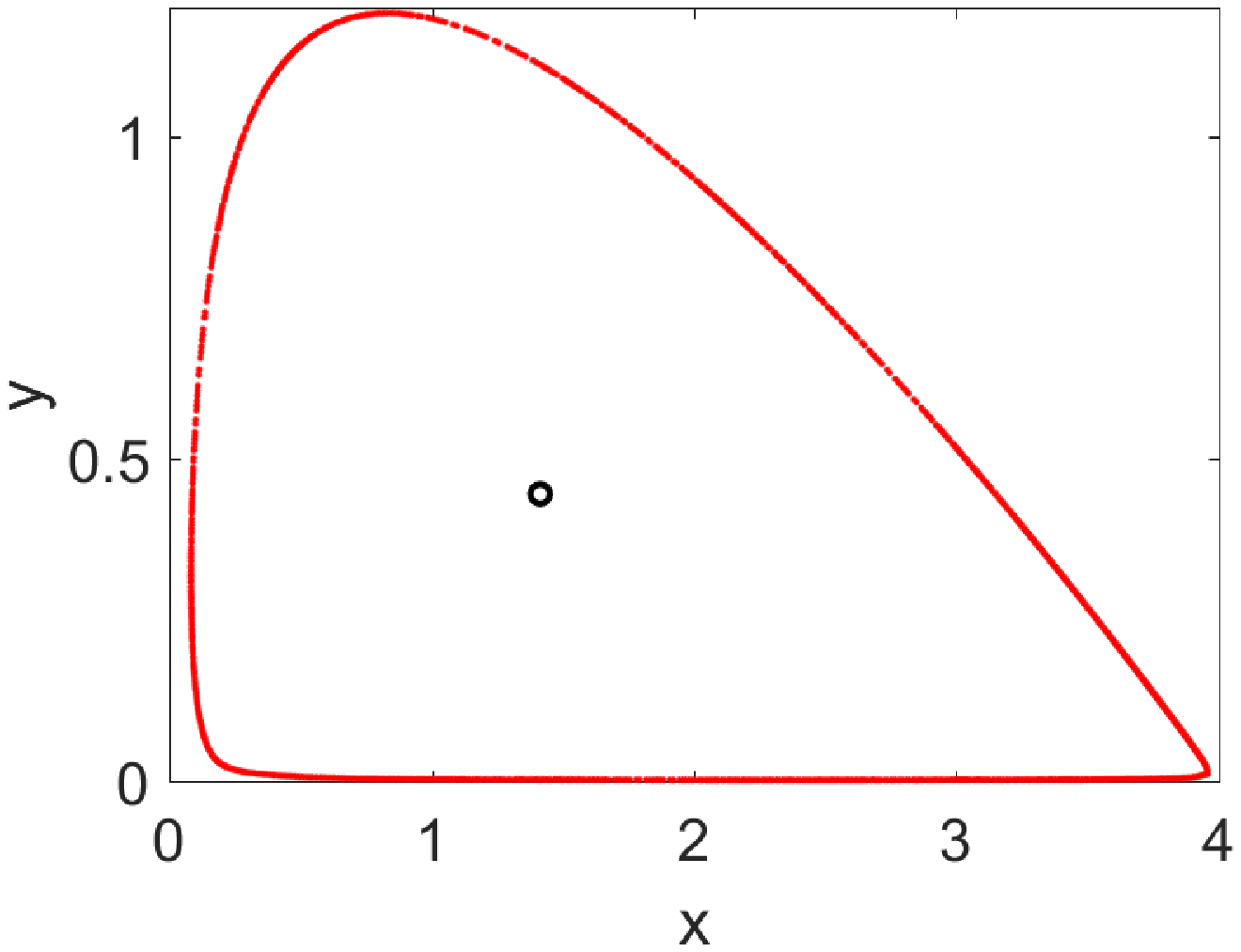} \\
(a) & (b)
\\
\includegraphics[height=2.0in, width=2.5in]{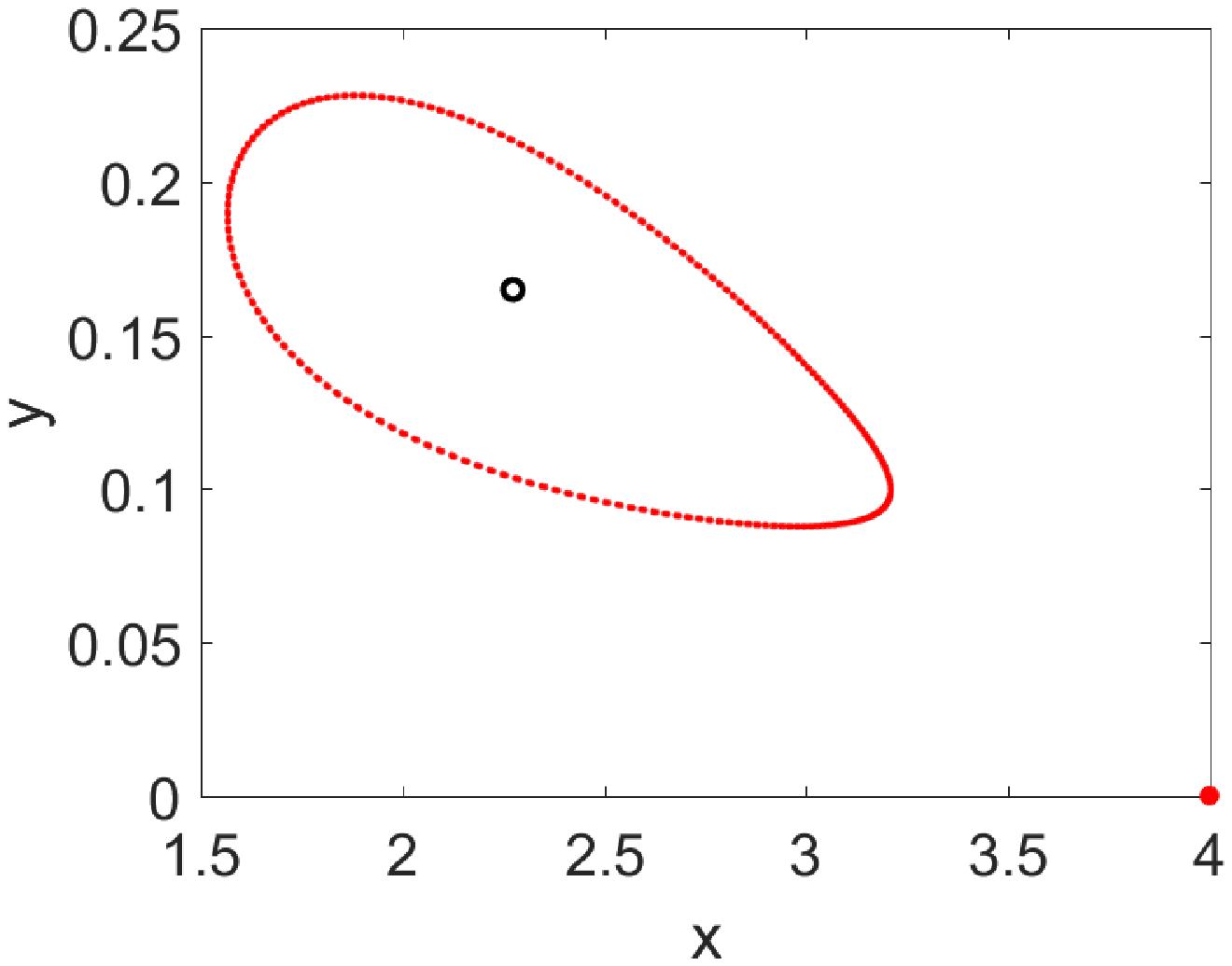} &
\includegraphics[height=2.0in, width=2.5in]{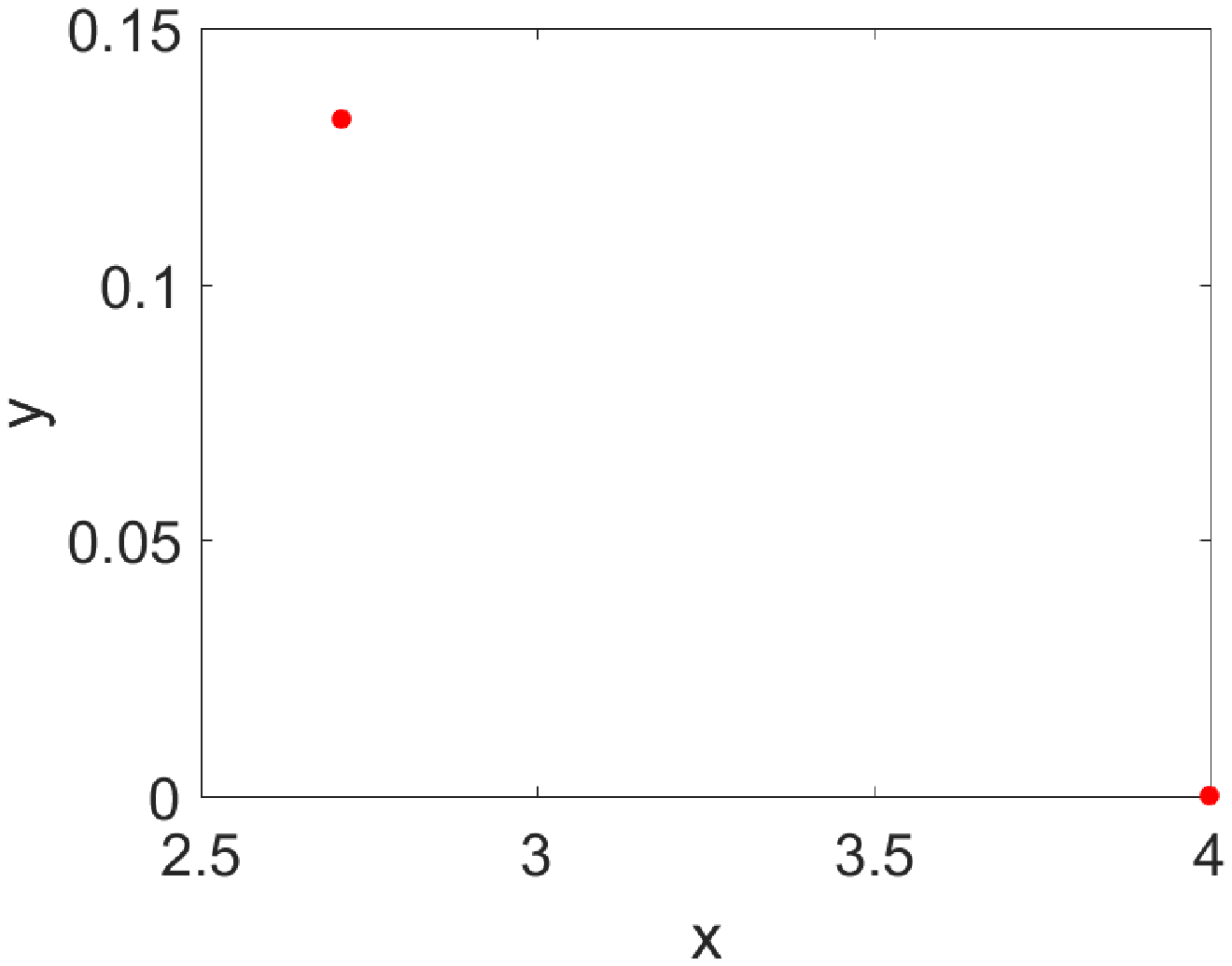}\\
(c) & (d)
\end{tabular}
\end{center}
\caption{Invariant closed curves and bistability are presented. In
(a), $2\alpha\leq 1$ and $\beta \bar x>1$. In (b), $2\alpha>1$ and
$\beta \bar x>1$. In (c) where $2\alpha>1$ and $\beta \bar x<1$,
the system has two interior steady states and two attractors are
shown using two different initial conditions. One solution
converges to the boundary steady state $E_1=(4,0)$ and the other
converges to the closed invariant circle. We decrease $\beta$ to
$\beta=3.76/20$ so that system \eqref{rmodel} still has two
interior steady states, where one is unstable and the other is
asymptotically stable. Using the same initial conditions as in
(c), one solution converges to $E_1$ while the other converges to
the interior steady state as shown in (d). }
\end{figure}
\end{document}